\documentclass[fleqn,journal]{ieeetran}
\usepackage{graphicx,amsmath,amssymb,arydshln,color,slashbox,cite}

\newcommand{\comment}[1]{}
\newtheorem{theorem}{Theorem}[section]
\newtheorem{corollary}[theorem]{Corollary}
\newtheorem{lemma}[theorem]{Lemma}

    \newcommand{\req}[1]{(\ref{#1.eq})}
    
    \newcommand{\be}{\begin{equation}}
    \newcommand{\ee}{\end{equation}}
    \newcommand{\beas}{\begin{eqnarray*}}
    \newcommand{\eeas}{\end{eqnarray*}}
    \newcommand{\bea}{\begin{eqnarray}}
    \newcommand{\eea}{\end{eqnarray}}
    
    \newcommand{\vt}{{\tilde{v}}}
    \newcommand{\xt}{{\tilde{x}}}
    
    \newcommand{\yh}{{\hat{y}}}
    \newcommand{\zh}{{\hat{z}}}

    \newcommand{\deltab}{{\boldsymbol{\delta}}}

    \newcommand{\ah}{\hat{a}}
    
    \newcommand{\ch}{\hat{c}}
    
    \newcommand{\ph}{\hat{p}}
    
    \newcommand{\fh}{\hat{f}}
    \newcommand{\gh}{\hat{g}}
    \newcommand{\gammah}{\hat{\gamma}}
    \newcommand{\phih}{\hat{\phi}}
    \newcommand{\Ah}{\hat{A}}
    \newcommand{\Bh}{\hat{B}}
    \newcommand{\Ch}{\hat{C}}
    \newcommand{\Hh}{\hat{H}}
    \newcommand{\Ph}{\hat{P}}
    \newcommand{\Kh}{\hat{K}}
    \newcommand{\Xh}{\hat{X}}
    \newcommand{\Zh}{\hat{Z}}
    \newcommand{\Yh}{\hat{Y}}
    \newcommand{\Fh}{\hat{F}}
    
    \newcommand{\Oh}{\hat{O}}
    
    \newcommand{\Gh}{\hat{G}}

    \newcommand{\Nb}{{\bar{N}}}
    
    \newcommand{\vbar}{\bar{v}}
    
    \newcommand{\cB}{{\cal B}}

    \newcommand{\bfo}{{\bf1}}
    \newcommand{\bfob}{{\bar{\bf1}}}
    \newcommand{\R}{{\mathbb R}}
    \newcommand{\C}{{\mathbb C}}
    \newcommand{\Z}{{\mathbb Z}}
    \newcommand{\expec}[1]{ {\cal E} \left\{ #1 \right\} }

\newcommand{\matbegin}{\left[ }
\newcommand{\matend}{\right] }
\newcommand{\obt}[2]{\matbegin \begin{array}{cc}
        	#1 & #2    \end{array} \matend }
\newcommand{\tbo}[2]{\matbegin \begin{array}{c}
        	#1 \\ #2    \end{array} \matend }
\newcommand{\tbt}[4]{\matbegin \begin{array}{cc}
        	#1 & #2 \\ #3 & #4  \end{array} \matend }
\newcommand{\thbo}[3]{
	  \matbegin \begin{array}{c}
       	#1 \\ #2 \\ #3
       	\end{array} \matend }
\newcommand{\obth}[3]{
  	\matbegin \begin{array}{ccc}
       	#1 & #2 & #3
       	\end{array} \matend }

\title{Coherence in Large-Scale Networks: \\
 Dimension-Dependent Limitations of Local Feedback}

\author{Bassam Bamieh, \IEEEmembership{Fellow, IEEE},
	Mihailo R.\ Jovanovi\'c, \IEEEmembership{Member, IEEE},
	Partha Mitra, and Stacy Patterson
\thanks{B. Bamieh is with the Department of Mechanical Engineering, University of California,
        Santa Barbara, CA 93106, USA
        ({\tt\small bamieh@engr.ucsb.edu}).}%
\thanks{M.\ R.\ Jovanovi\'c is with the Department of Electrical and Computer Engineering,
        University of Minnesota, Minneapolis, MN 55455, USA
        ({\tt\small mihailo@umn.edu}).}%
\thanks{Partha Mitra is with Cold Spring Harbor Laboratory,
        Cold Spring Harbor, NY, USA
        ({\tt\small mitra@cshl.edu}).}%
\thanks{S. Patterson is with the Department of Electrical Engineering, Technion - Israel Institute of Technology,
Haifa, 32000, Israel
        ({\tt\small stacyp@ee.technion.ac.il}).}%

     \thanks{This research is partially supported by NSF grants ECCS-0802008 and CMMI-0626170, and AFOSR FA9550-10-1-0143.}
}

\date{Feb 17, 2009}

\begin{document}

\maketitle

\begin{abstract}
We consider distributed consensus and vehicular formation control problems. Specifically we address the question of whether local feedback is sufficient to
maintain coherence in large-scale networks subject to stochastic disturbances. We define macroscopic performance measures which are global quantities that capture the notion of coherence; a notion of global order that
quantifies how closely the formation resembles a solid object.
We consider how these measures scale asymptotically with network size in the topologies of regular lattices in 1, 2 and higher dimensions, with vehicular platoons corresponding to the
1 dimensional case.
A common
phenomenon appears where a higher spatial dimension implies a more favorable
scaling of coherence measures, with a dimensions of 3 being necessary to achieve coherence in consensus and vehicular formations under certain conditions. In particular, we show that it is
impossible to have large coherent one dimensional vehicular platoons with only local feedback.
We analyze these effects in terms of the underlying energetic modes of motion, showing that they take the form of large temporal and spatial scales resulting in an accordion-like motion of formations. A conclusion can be drawn that in low spatial dimensions, local feedback is unable to regulate large-scale disturbances, but it can in higher spatial dimensions. This phenomenon is distinct from, and unrelated to string instability issues which are commonly encountered in
control problems for automated highways.
\end{abstract}

\section{Introduction}
The control problem for strings of vehicles (the so-called platooning problem) has been
extensively studied in the last two decades, with original problem formulations and studies dating
back to the
60's~\cite{shldeshedtomwalzhamcmpenshemck91,swahed99,swahed96,melkuo71,levath66}. These problems are also intimately related to more recent formation
flying and formation control problems~\cite{marcorbull07}.
It has long been observed in platooning
problems that to achieve reasonable performance, certain global information such as leader's
position or state need to be broadcast to the entire formation. A precise analysis of the limits
of performance associated with localized versus global control strategies does not appear
to exist in the formation control literature. In this paper we study the platooning problem as
the 1 dimensional version of a more general vehicular formations control problem on regular lattices
in arbitrary spatial dimensions.
For such problems,
we investigate the limits of performance of any local feedback law that is globally stabilizing.
In particular, we propose and study measures of the  coherence of the
formation. These are measures that capture the notion of how well the formation resembles
a rigid lattice or a solid object.

The coherence of a formation is a different concept from, and often unrelated to, string instability.
In the platooning case (i.e.\ 1 dimensional formations), which turns out to be most problematic, a localized
feedback control law may posses string stability in the sense that the effects of any injected
disturbance do not grow with spatial location. However, as we show in this paper, it is impossible
to achieve a large coherent
formation with only localized feedback if all vehicles are subject to any amount
of distributed stochastic  disturbances. The net effect is that with the best localized feedback, a 1 dimensional
formation will appear to behave well on a ``microscopic'' scale in the sense that distances between
neighboring vehicles will be well regulated. However, if a large formation is observed in its entirety,
it will appear to have temporally slow, long spatial wavelength modes that are unregulated,
resembling an ``accordion''
type of motion. This is not a safety issue, since the formation is microscopically well regulated,
but it might effect throughput performance in a platooning arrangement since throughput does depend
on the coherence or rigidity of the formation.

The phenomenon that we discuss occurs in both consensus algorithms and vehicular
formation problems. We therefore treat both as instances of networked dynamical systems with
first order and second order local dynamics respectively.
Both problems are set up in the $d$-dimensional torus $\Z_N^d$.
We begin in section~\ref{probform.sec}
with problem formulations of the consensus type and vehicular formations,
where we view the former as a first order dynamics version of the latter. In section~\ref{perf.sec}, we
define macroscopic and microscopic measures of performance in terms of variances of various
quantities across the network. We argue that the macroscopic measures capture the notion
of coherence. We also present compact formulae for calculating those measures as $H^2$ norms of
systems with suitably defined output signals. These norms are calculated using traces of
system Grammians, which in turn are  related to sums involving eigenvalues of the underlying system
and feedback gains matrices. Since the network topologies we consider are built over Tori networks, these
system matrices are multi-dimensional circulant operators, and  their eigenvalues are calculated
as the values of the Fourier symbols of the underlying feedback operators, thus allowing for a rather
direct relation between the structure of the feedback gains and the system's norms.
Much of the remainder of the paper is devoted to establishing
asymptotic (in network size) bounds
for these performance measures for each underlying spatial dimension. Section~\ref{upbounds.sec}
establishes
upper bounds of standard algorithms, while section~\ref{lbounds.sec}
is devoted to establishing lower bounds
for {\em any} algorithm
that satisfies a certain number of structural assumptions including the locality
of feedback and boundedness of control effort.
This shows that asymptotic limits of performance are determined by the network
structure rather than the selection of parameters of the feedback algorithm.
We pay particular attention to the role of control effort as our
lower bounds are established for control laws that have bounded control effort in a stochastic sense.
Some numerical examples illustrating the
lack-of-coherence phenomenon are presented in section~\ref{comp.sec},
as well as an illustration of how it is distinct from
string instability. The interested reader may initially skim this section which numerically illustrates the
basic phenomenon we study analytically in the remainder of the paper.
We end in section~\ref{disc.sec} with a discussion of related work in which various versions of this phenomenon were observed, as well as a discussion of some open questions.

\subsection*{Notation and Preliminaries}

The networks we consider are built over the d-dimensional Torus $\Z_N^d$. The one-dimensional
Torus $\Z_N$ is simply the set of integers $\{0,1,\ldots,N-1\}$ with addition {\em modulo N}
($\mod N$), and $\Z_N^d$ is
the direct product of $d$ copies of $\Z_N$. Functions defined on $\Z_N^d$ are called {\em arrays},
and we use multi-index notation for them, as in $a_k ~=~ a_{(k_1,\ldots,k_d)}$ to denote
individual entries of an array. Indices are added in the $\Z_N^d$ arithmetic as follows
\begin{eqnarray*}
	(r_1,\dots,r_d) & = & (k_1,\dots,k_d)  ~+~(l_1,\dots,l_d) 		\\
	& \Updownarrow & 	\\
	r_i & = &  \left(k_i ~+~ l_i \right)_N , ~~~~~i=1,\ldots,d,
\end{eqnarray*}
where $()_N$ is the operation $\mod N$.
The set $\Z_N^d$ and the corresponding addition operation can be
visualized as a ``circulant'' graph in d-dimensional space with edge nodes connected to nodes on corresponding
opposite edge of the graph.

The multi-dimensional Discrete Fourier Transform is used throughout. All states are
multi-dimensional arrays which we define as real or complex vector-valued functions on
the Torus $\Z_N^d$. The Fourier transform (Discrete Fourier Transform) of an array $a$ is denoted with $\ah$.
We refer to indices of spatial Fourier transforms
as {\em wavenumbers}. Generally, we use $k$ and $l$ for spatial indices and
$n$ and $m$ for wavenumbers. For example, an array $a_{(k_1,\dots,k_d)}$ has $(k_1,\dots,k_d)$ as
the spatial index, while its Fourier transform $\ah_{(n_1,\ldots,n_d)}$ has the index $(n_1,\dots,n_d)$
as the wavenumber. The wavenumber is simply a spatial frequency variable.
Some elementary properties of this Fourier transform are
summarized in Appendix~\ref{fourier.appen},

Convolution operators arise naturally over $\Z_N^d$. Let $a$ be any array of numbers (or matrices)
over $\Z_N^d$,
that is $a:\Z_N^d \rightarrow \C $ (or $\C^{n\times n}$).
Then the operator
 $T_a$ of multi-dimensional circular convolution with the array $a$ is defined as follows
 \begin{eqnarray*}
 	g & = &  T_a f ~=~ a \star f		\\
	& \Updownarrow &  			\\
	g_{(k_1,\ldots,k_d)} & = &  \sum_{(l_1,\ldots,l_d)\in \Z_N^d}
	a_{(k_1,\ldots,k_d)-(l_1,\ldots,l_d)} ~f_{(l_1,\ldots,l_d)}.
\end{eqnarray*}
Note that $f$ and $g$  may be scalar or vector-valued (depending on whether $a$ is scalar or
matrix-valued respectively), and that the arithmetic for $(k_1,\ldots,k_d)-(l_1,\ldots,l_d)$ is done
in $\Z_N^d$, i.e. arithmetic $\mod N$ in each index as described above.

 It is important to distinguish between an array $a$ and the corresponding linear operator $T_a$.
The Fourier
transform $\ah$ of the array $a$ is called the {\em Fourier symbol} of the operator $T_a$. It is a standard
fact that the eigenvalues of the operator $T_a$ are exactly the values of the Fourier transform $\ah$,
i.e. the values of its Fourier symbol. When $a$ is matrix valued, then the eignvalues of $T_a$ are the union
of all eigenvalues of $\ah_{(n_1,\ldots,n_d)}$ as the wavenumber $(n_1,\ldots,n_d)$ runs through $\Z_N^d$,
i.e.
\[
	\sigma\left(T_a \right) ~=~ \bigcup_{(n_1,\ldots,n_d)\in \Z_N^d} \sigma \left( \ah_{(n_1,\ldots,n_d)} \right),
\]
where $\sigma(.)$ refers to the spectrum of a matrix or operator (all finite-dimensional in our case).

In this paper, we use the term {\em dimension} to refer exclusively to the spatial dimension of underlying
networks. To avoid confusion with the notion of state dimension, we refer to the dimension of the state space
of any dynamical system as the {\em order} of that dynamical system.

The vector dimension of signals is mostly suppressed to keep the notation from being cumbersome. For
example, the state of node $(k_1,\ldots,k_d)$ in the d-dimensional Torus is written as
	\[
		x_{(k_1,\ldots,k_d)}(t).
	\]
It is a scalar-valued signal for consensus problems, and vector-valued (in $\R^d$) signal for vehicular
formation problems.

We use $M^T$ to denote the transpose of a matrix $M$, and $M^*$ to denote the complex-conjugate transpose
of a matrix $M$ or the adjoint of an operator $M$. Although all operators in this paper are finite dimensional,
we sometimes refer to them as operators rather than matrices since we often avoid writing
the cumbersome explicit
matrix representations (such as in the case of multi-dimensional convolution operators).

\section{Problem Formulation}
	\label{probform.sec}

We formulate two types of problems, consensus and vehicular formations. The mathematical setting is analogous in both problems, with the main difference being that vehicular models have two states (position and velocity) locally at each site in contrast to a scalar local state in consensus problem. This difference
leads to more severe asymptotic scalings in vehicular formations as will be shown in the sequel.

\subsection{Consensus with random insertions/deletions}

We begin by formulating a continuous-time version of the  consensus algorithm
with additive stochastic disturbances in the
dynamics~\cite{X07,patbamela}.
As opposed to standard consensus algorithms without additive disturbances, nodes do not achieve
equilibrium asymptotically but fluctuate around the equilibrium, and the variance of this
fluctuation is a measure of how well approximate consensus is achieved.
This formulation can be used to model scenarios such as
load balancing over a distributed file system, where the additive noise represents file insertion and deletion,
parallel processing systems where the noise processes model job creation and completion,
or flocking problems in the presence of random forcing disturbances.

We consider a consensus algorithm over undirected tori, $\Z_N^d$, where
the derivative of the scalar state at each node is determined as a weighted average of the differences between that node and all its $2d$ neighbors.
One possible such algorithm is given by
\bea
	\dot{x}_k
	& = &  \beta \left[ \left(  x_{(k_1-1,\ldots,k_d)} - x_k \right)			  +\cdots +  \right.  	 \nonumber	\\
	& &    \left. ~~  \left( x_{(k_1,\ldots,k_d+1)}  - x_k \right)
						 \right]
		~+~ w_k, 						\nonumber			\\
	& = &    (-2d\beta) ~x_k +	\beta \left[
		x_{(k_1-1,\ldots,k_d)} ~+~ x_{(k_1+1,\ldots,k_d)}
		+		\right.	\nonumber	\\
	& &
		\left.		\cdots + x_{(k_1,\ldots,k_d-1)} ~+~ x_{(k_1,\ldots,k_d+1)}
						\right]
		~+~ w_k, 											\label{cncsform.eq}
\eea
where we have used equal weights $\beta > 0$ for all the differences. The process disturbance $w$ is a mutually uncorrelated white stochastic process.
We call this the {\em standard consensus algorithm\/} in this paper
since it is essentially the same as other well-studied consensus
algorithms~\cite{T84, B90,jadlinmor03,PBE06, carli2007cca}.

The sum in the equation above can be written as a multimulti-dimensional convolution by defining
the array
\be
	O_{(k_1,\ldots,k_d)} = \left\{
		\begin{array}{ll}
			-2d\beta &  k_1=\cdots=k_d=0,	\\
			\beta & k_i =\pm 1, ~\mbox{and}~ k_j=0 ~\mbox{for}~ i\neq j, 	\\
			0 & \mbox{otherwise}.
		\end{array}		\right.
   \label{scp_array.eq}
\ee
The system~\req{cncsform} can then be written as
\be
	\dot{x} ~=~ O\star x ~+~ w,
   \label{stand_cncs.eq}
\ee
where $\star$ is circular convolution in $\Z_N^d$.

We recall that we use the
operator notation $T_a ~x ~:=~ a\star x$ to indicate the circulant operator of
convolution with any array $a$.  With this notation, a general spatially invariant consensus
algorithm can be written
abstractly as
\be
	\dot{x} ~=~ T_a ~x ~+~ w,
   \label{gen_cncs.eq}
\ee
for any array $a$ defined over $\Z_N^d$.
Such algorithms
can be regarded as a combination of open loop dynamics
\[
	\dot{x}_k ~=~ u_k ~+~ w_k , ~~~~~~k\in\Z_N^d,
\]
with the feedback ``control'' $u ~=~ T_a x$, where the feedback operator array is to be
suitably designed. With this point of view, consensus algorithms can be thought of as first
order dynamics versions of vehicular formation problems that we introduce next.

\subsection{Vehicular Formations}

Consider $N^d$ identical vehicles arranged in a $d$-dimensional torus, $\Z_N^d$, with the
double integrator dynamics
\be
	\ddot{x}_{(k_1,\ldots,k_d)} ~=~ u_{(k_1,\ldots,k_d)} ~+~ w_{(k_1,\ldots,k_d)},
   \label{xeqmo.eq}
\ee
where $(k_1,\ldots,k_d)$ is a multi-index with each $k_i\in\Z_N$, $u$ is the control
input and $w$ is a mutually uncorrelated white stochastic process which can be considered to model random forcing. In the sequel, we will also consider the consequences of the presence of viscous
friction terms in models of the form
\be
	\ddot{x}_{(k_1,\ldots,k_d)} ~=~ -\mu \dot{x}_{(k_1,\ldots,k_d)} ~+~
		u_{(k_1,\ldots,k_d)} ~+~ w_{(k_1,\ldots,k_d)},
   \label{wviscfric.eq}
\ee
where $\mu>0$ is the linearized drag coefficient per unit mass.

Each position vector
$x_k$ is a $d$-dimensional vector with components $x_k ~=~\obth{x_k^1}{\cdots}{x_k^d}^T$.
The objective is to have the $k$th vehicle in the formation follow the desired trajectory
$\bar{x}_k$
\[
	\bar{x}_k ~:=~ \vbar t  ~+~ k\Delta   ~~\Leftrightarrow~~
	\thbo{\bar{x}_k^1}{\vdots}{\bar{x}_k^d} ~:=~ \thbo{\vbar^1}{\vdots}{\vbar^d} t
		~+~ \thbo{k_1}{\vdots}{k_d}\Delta,
\]
which means that all vehicles are to move with constant heading velocity $\vbar$ while maintaining
their respective position in a $\Z_N^d$ grid with spacing of $\Delta$ in each dimension. The situation
of different spacings in different directions can be similarly represented, but is not considered for notational
simplicity.

The deviations from desired trajectory are defined as
\[
	\tilde{x}_k ~:=~ x_k ~-~\bar{x}_k, ~~~~~\vt_k ~:=~ \dot{x}_k - \vbar.
\]
We assume the control input to be full state feedback and linear in the variables $\xt$
and $\vt$ (therefore affine linear in $x$ and $v$), i.e. $u ~=~ G \xt ~+~F \vt$, where $G$
and $F$ are the linear feedback operators.
The equations of motion for the controlled system are thus
	\be
		\frac{d}{dt} \tbo{\xt }{\vt} ~=~ \tbt{0}{I}{G}{F} \tbo{\xt}{\vt} ~+~\tbo{0}{I} w.
	   \label{vfsys.eq}
	\ee
We note that the above equations are written in operator form, i.e. by suppressing the
spatial index of all the variables.

\textit{Example:}
The operators $G$ and $F$ will have some very special structure depending on assumptions
of the type of feedback and measurements available. Consider for example a feedback control
of the $k$th vehicle (in a one dimensional formation) of the following form
	\[
	\begin{array}{rclcll}
		u_k & = &  {g_{_+} (x_{k+1}-x_k-\Delta)}  & + &
				 {g_{_-} (x_{k-1} - x_k -\Delta)}  & + \\
			& & 	   {f_{_+} (v_{k+1}-v_k)}          & + &
					{f_{_-}(v_{k-1} - v_k)}  & +  \\
			& & 		{ g_o \left(x_k - \bar{x}_k \right)}  & +
            & { f_o \left(v_k - \vbar \right)}, &
	\end{array}
	\]
where the  $g$'s and $f$'s  are design constants.
The first two lines represent look-ahead and look-behind position and velocity error feedbacks
respectively. We refer
to such terms as {\em relative feedback} since they {\em only\/} involve measurements of differences.
On the other hand, terms in the last line require knowledge of positions and velocities in an absolute
coordinate system (a grid moving at constant velocity),
and we thus refer to such terms as {\em absolute feedback}.
For later reference, it is instructive to write the feedback in the above example in terms of the state variables $\xt$ and $\vt$ as
	\begin{eqnarray}
		u_k & = &  {g_{_+} (\xt_{k+1}-\xt_k)}   ~+~
				 {g_{_-} (\xt_{k-1} - \xt_k)}   ~+ 			\nonumber \\
			& & 	   {f_{_+} (\vt_{k+1}-\vt_k)}         ~~  + ~
					{f_{_-}(\vt_{k-1} - \vt_{k})}   ~+  		\nonumber \\
			& & 		{ g_o ~\xt_k}   ~+~
			 		{ f_o~  \vt_k }.
			  	\label{platoonfdbk.eq}
	\end{eqnarray}

\subsection{Structural assumptions}

We now list the various assumptions that can be imposed on system operators
and on the control feedbacks $G$ and $F$. These are structural restrictions representing
the structure of open loop dynamics and measurements, and the type of feedback control
available respectively.
						\newcounter{Acount}
\begin{list}
		{	\bf (A\arabic{Acount})	 }
		{ 	\usecounter{Acount}
			\setlength{\itemindent}{-0.3em}
			\setlength{\leftmargin}{1em}		}
	\item \textbf{Spatial Invariance.}
		All operators
		are spatially invariant with respect to $\Z_N^d$. This
		implies that they are convolution operators. For instance, the operation $Gx$
		can be written as the convolution (over $\Z_N^d$) of the array $x$ with
		an array $\{G_{(k_1,\ldots,k_d)} \}$
		\be
			(Gx)_k ~=~ \sum_{l\in\Z_N^d} G_{k-l} ~x_l,
		   \label{mvecconv.eq}
		\ee
		where the arithmetic for $k-l$ is done in $\Z_N^d$. For each $k$, the array element $G_k$
		is a $d\times d$ matrix ($G$ is then an $N^d \times N^d$ operator).
		Note that in the absence of spatial invariance, each term of the sum
		in~\req{mvecconv} would need to be written as $G_{k,l} ~x_l$. That is, one would require a
		two-indexed array of matrices $G_{k,l}$ rather than a single-indexed array.
		
In the example above of a one dimensional circular formation,
the array elements for position feedback are given by
$\left\{ \left( g_o-g_{_+}-g_{_-} \right)   ,  g_{_+}, 0, \ldots, 0, g_{_-} \right\}$.
	
\item \textbf{Relative vs.\ Absolute Feedback.}		
		We use the term Relative Feedback when given feedback involves only differences between
		quantities.
		For example, in position feedback, this implies that
		for each term of the form $\alpha x_{(k_1,\ldots,k_d)}$ in the convolution,
		another term of the form $-\alpha x_{(l_1,\ldots,l_d)}$ occurs for some other
		multi-index $l$. This implies that the array $G$ has the property
		\be
			\sum_{k\in\Z_N^d} G_k ~=~  0.
		   \label{relmeas.eq}
		\ee
		We use the term Absolute Feedback when given operator does not satisfy this assumption.
		
		Note that in the example above, relative position feedback corresponds to $g_o=0$, and in this case,
		condition~\req{relmeas} is satisfied.
	
\item \textbf{Locality.}
		The feedbacks use only local information from a neighborhood of width $2q$,
		where $q$ is independent of $N$. Specifically,
		\begin{eqnarray}
			G_{(k_1,\ldots,k_d)} & = &  0, 	~~\mbox{if} ~k_i> q, ~\mbox{and} ~k_i <N-q  \nonumber \\
				& &  \mbox{for any $i\in\{1,\ldots,d\}$}. ~
			\label{locality.eq}
		\end{eqnarray}
		The same condition holds for $F$.
	
\item \textbf{Reflection Symmetry.}
		The interactions between vehicles have mirror symmetry.
		This has the consequence that the arrays
		representing $G$ and $F$ have even symmetry, e.g. for each nonzero
		term like $\alpha G_{(k_1, \ldots,k_d)}$ in the array there is a corresponding
		term $\alpha G_{(-k_1, \ldots,-k_d)}$. This in particular implies that the Fourier
		symbols of $G$ and $F$ are real valued.
		
		In the example above, this condition gives $g_{_+}=g_{_-}$ and $f_{_+}=f_{_-}$.
	
\item \textbf{Coordinate Decoupling.} For $d\geq 2$, feedback control of thrust in each
		coordinate direction depends only on measurements of position and velocity error
		vector components in
		that coordinate. This is equivalent to imposing that each array element $G_k$ and $F_k$
		are $d\times d$ diagonal matrices. For further
		simplicity we assume those diagonal elements to be
		equal, i.e.
		\be
			G_k ~=~ \textrm{diag}\{g_k, \ldots,g_k\}, ~~
			F_k ~=~ \textrm{diag}\{f_k, \ldots,f_k\}.
		   \label{diagfdbk.eq}
		\ee
		This in effect renders the matrix-vector convolution in~\req{mvecconv} into $d$ decoupled scalar
		convolutions.
\end{list}	

Assumptions {\bf (A1)} through {\bf (A3)} appear to be important for subsequent developments, while assumptions {\bf (A4)} and {\bf (A5)} are made to simplify calculations.

\section{Performance Measures}
	\label{perf.sec}

We will consider how various performance measures scale with system size for the consensus
and vehicle formations problems. Some of these measures can be quantified as
steady state variances of outputs of linear systems driven by stochastic inputs, so we
consider some generalities first.
Consider a general linear system driven by zero mean white noise with unit covariance
\begin{eqnarray*}
	\dot{x} & = &  Ax ~+~ Bw,	\\
	y & = & Hx.
\end{eqnarray*}
Since we are interested in cases where $A$ is not necessarily Hurwitz (typically due to a single
unstable mode at the origin representing motion of the mean), the state
$x$ may not have finite steady state variances. However, in all cases we consider
here the outputs $y$ do have finite variances, i.e. the unstable modes of $A$
are not observable from $y$. In such cases, the output does have a finite steady state variance,
which is quantified by the square of the $H^2$ norm of the system from $w$ to $y$
    \be
	V
	~ := ~ \sum_{k\in \Z_N^d} \lim_{t\rightarrow\infty} E\left\{y^*_k(t) y_k(t)\right\}  ,
   \label{Vdef.eq}
\ee
where the index $k$ ranges over all ``sites'' in the lattice $\Z_N^d$.

We are interested in spatially invariant problems over discrete Tori.
This type of invariance implies that the variances
of all outputs are equal, i.e. $E\left\{y^*_k y_k\right\}$ is independent of $k$.  Thus,
if the output variance at a given site is to be computed, it is simply
the total $H^2$ norm divided by the system size
\be
	E\left\{y_k^*y_k\right\} ~=~
		\frac{1}{M}\sum_{l\in\Z_N^d} E\left\{y_l^*y_l\right\} ~=~ \frac{V}{M},
   \label{iov.eq}
\ee
where $M$ is the size of the system ($M = N^d$ for $d$-dimensional Tori). We refer to quantities
like~\req{iov} as {\em individual output variances}.

Next, we define several different performance measures and give the corresponding
output operators for each measure for both the consensus and vehicular formation problems. In the vehicular formation problem, we assume for simplicity that the output involves positions only,
 and thus the output equation has the form
\[
y =  \left[\begin{array}{cc}C & 0\end{array}\right] \tbo{\xt}{\vt},
\]
i.e. $H=\obt{C}{0}$,
where $C$ is a circulant operator. A consensus problem with the same performance measure
has a corresponding output equation of the form (with the same $C$ operator)
\[
y = C x.
\]

\textit{Performance Measures:}
We now list the three different performance measures we consider.
						\newcounter{Pcount}
\begin{list}
		{	{\bf (P\arabic{Pcount})}	 }
		{ 	\usecounter{Pcount}
			\setlength{\itemindent}{-0.3em}
			\setlength{\leftmargin}{1em}		}

\item \textbf{Local error.}  This is a measure of the difference between neighboring
	nodes or vehicles. For the consensus problem, the $k$th output (in the case of one dimension) is
	defined by
	\[
	y_k := x_k  - x_{k-1}.
	\]
	For the case of vehicular formations,  local error is
	the difference of neighboring vehicles positions from
	desired spacing, which can equivalently be written as
	\[
	y_k := \xt_k - \xt_{k-1}.
	\]
	The output operator is then given by $C~:=~ (I-D)$, where $D$ is the right shift
	operator,  $(Dx)_k ~:=~ x_{k-1}$.

	In the case of $d$ dimensions, we define a vector output that contains as components the local error
	in each respective dimension, i.e.
	\be
		C ~:=~ \frac{1}{\sqrt{2d}}  \obth{I~-D^1}{\cdots}{I~-~D^d}^T,
	   \label{Clocaldef.eq}
	\ee
	where $D^r$ is the right shift along the $r$th dimension, i.e.
	$(D^rx)_{(k_1,\ldots,k_r,\ldots,  k_d)}  ~:=~ x_{(k_1,\ldots,k_r-1,\ldots,  k_d)}  $, and $1/\sqrt{2d}$ is a
	convenient normalization factor. This operator is closely related to the standard consensus operator $O$ in Eq.~\req{scp_array} by the following easily established identity
	\be
		C^*C ~=~ \frac{-1}{2d\beta} ~O.
	   \label{Oiden.eq}
	\ee
	
\item \textbf{Long range deviation (Disorder).} In the consensus problem, this corresponds to measuring
	the disagreement between the two furthest nodes in the network graph.
	Assume for simplicity that $N$ is even and we are in dimension 1.
	Then, the most distant node from node $k$
	is $\frac{N}{2}$ hops away, and we define	
	\[
	y_k := x_k - x_{k+\frac{N}{2}}.
	\]
	
	In the vehicular formation problem,
	long range deviation corresponds to measuring the deviation
	of the distance between the two most distant vehicles from what it should be.
	The most distant vehicle to the $k$th one is the vehicle indexed by $k+\frac{N}{2}$. The desired distance between them is $\Delta \frac{N}{2}$,
	and the deviation from this distance is
	\be
		y_k := x_k-x_{k+\frac{N}{2}}-\Delta \frac{N}{2} ~=~ \xt_k - \xt_{k+\frac{N}{2}}.
	   \label{lrd.eq}
	\ee
	We consider the variance of this quantity to be a measure of
	disorder, reflecting the lack of ``end-to-end rigidity'' in the vehicle formation.
	
	Generalizing this measure to $d$
	dimensions yields an output operator of the form
	\be
		C ~:=~ T_{(\deltab^0-\deltab^{(N/2, \ldots, N/2)})},
	    \label{Clrd.eq}
	\ee
	i.e. the operator of convolution with the
	array\footnote{By a slight abuse of notation, we define the shifted Kronecker delta
	$\deltab^l_{k} := \delta_{k-l}$, where $\delta_{k}=1$ for $k=0$, and zero otherwise, is the
	standard Kronecker delta.  With this notation, $\deltab^0$ is also the standard Kronecker delta.}
	  $\deltab^0-\deltab^{(N/2, \ldots, N/2)}$.
	\comment{
	
	The Fourier
	symbol of $C_1$ is then easily calculated to be
	\beas
		\Ch_1
		& = &  1 -  e^{-i\frac{2\pi}{N}(\frac{N}{2}n_1+\cdots+\frac{N}{2}n_d)}      \\
		 & = &   1 -  e^{-i{\pi}(n_1+\cdots+n_d)}     \\
		& = &  \left\{ \begin{array}{lcl}
				0, & ~& (n_1+\cdots+n_d) ~~\mbox{even} \\
				2, & & (n_1+\cdots+n_d) ~~\mbox{odd} \\
				\end{array}
			\right. .
	\eeas
	Putting this together with~\req{Vv} gives the $H^2$ norm for long range deviation
	(which we will refer to as $V_{lrd}$)
	\be
		 V_{lrd} ~=~ 	
		 	\sum_{n_1+\cdots+n_d ~\mbox{odd},~ n\in\Z_N^d} \frac{1}{\Kh_n\Fh_n}.
   	   \label{Vv2.eq}
	\ee	
	}
\item \textbf{Deviation from average.}
	For the consensus problem, this quantity measures the deviation of each state
	from the average of all states,
	\be
		y_k := x_k - \frac{1}{M} \sum_{l \in\Z_N^d} x_l.
	   \label{davcons.eq}
	\ee
	In operator form we have $y =~ (I-T_\bfob)x$,
	where $\bfob$ is the array of all elements equal to $1/M$.

	In  vehicular formations, this measure can be interpreted as the deviation of each vehicle's position
		error from the average of the
		overall position error
		$y~=~ (I-T_\bfob) \xt$.
\end{list}

We note that performance measures {\bf (P1)} through {\bf (P3)} are such that $C$ can be represented as a convolution with an array $\{C_k\}$ which has the
property $\sum_{~ k\in\Z_N^d} C_k = 0$.
This condition causes the mean mode at zero to be unobservable, and thus guarantees that all outputs
defined above have finite variances.

We refer to the performance measure {\bf (P1)} as a {\em microscopic error} since it involves quantities local to any given site. This is in contrast to the measures {\bf (P2)} and {\bf (P3)} which involve quantities that are far apart
in the network, and we thus refer to these as {\em macroscopic errors}. We consider the macroscopic errors as measures of disorder or equivalently, lack of coherence. As we will show in the sequel, both macroscopic measures scale similarly asymptotically with system size, which justifies using either of them as a measure of disorder.

\subsubsection*{Formulae for variances}

Since we consider spatially invariant systems and in particular systems on the discrete Tori
$\Z_N^d$, it is possible to derive formulae for the above defined measures in terms of the
Fourier symbols of the operators $K$, $F$ and $C$. Recall the state space formula for the $H^2$
norm $V$ defined in~\req{Vdef}
\[
	V~ = ~ 	\mbox{tr}\left(
		\int_{0}^{\infty} B^*  e^{A^*t}H^*H e^{A t}  B ~dt
		\right).
\]

When $A$, $B$ and $H$ are circulant operators, traces can be rewritten
in terms
of their respective Fourier symbols (see ~\req{circtr})  as
\bea
	V
	& = &  \mbox{tr}\left( 	\sum_n
	\int_{0}^{\infty} \Bh^*_n e^{\Ah_n^*t}\Hh_n^*\Hh_n e^{\Ah_n t}  	\Bh_n ~d t       \right)
			\\
	& = &  \sum_n \mbox{tr}\left( \Bh^*_n \Ph_n \Bh_n \right) ,
   \label{totalh2.eq}
\eea
where the individual integrals are defined as
\begin{equation}
	\Ph_n 	~:=~ \int_0^\infty e^{\Ah_n^*t}\Hh_n^*\Hh_n e^{\Ah_n t} dt.
   \label{pkinteg.eq}
\end{equation}
If $\Ah_n$ is Hurwitz, then  $\Ph_n$ can be obtained by solving the Lyapunov equation
\begin{equation}
	\Ah_n^*\Ph_n+\Ph_n\Ah_n ~=~ -\Hh_n^*\Hh_n.
   \label{pklyap.eq}
\end{equation}
For wavenumbers $n$ for which $\Ah_n$ is not Hurwitz, $\Ph_n$ is still finite if  the non-Hurwitz
modes of $\Ah_n$ are not observable from $\Hh_n$. In this case
we can analyze the integral in~\req{pkinteg} on a case by case basis.

The Lyapunov equations are easy to solve in the Fourier domain. Equation~\req{pklyap} is a scalar
equation in the Consensus case and a $2d\times 2d$ matrix equation in the Vehicular
case\footnote{Note that in $d$ dimensions, the transformed state vector is of dimension
$2d$ for each wavenumber $n$.}.
The two respective
calculations are summarized in the next lemma. The proof is given in the Appendix.

    \begin{lemma}
The output variances~\req{Vdef} for the consensus and vehicular problems satisfying assumptions {\bf (A1)}-{\bf (A5)} are given by
	\label{V.lemma}
\be
	V_c ~=~ -\frac{1}{2} \sum_{n\neq 0,~ n\in\Z_N^d} \frac{|\ch_n|^2}{\Re(\ah_n)},
   \label{Vc.eq}
\ee
\be
	V_v ~=~ \frac{d}{2} \sum_{n\neq 0,~ n\in\Z_N^d} \frac{|\ch_{n}|^2}{\gh_n \fh_n  },
   \label{Vv.eq}
\ee
where $\Re(\ah_n)$ is the real part of $\ah_n$,
$\ch$ is the Fourier symbol of the output operator corresponding to the
performance index under consideration, and $\ah$, $\gh$ and $\fh$ are the Fourier symbols
of the consensus operator~\req{gen_cncs},
and the position and velocity feedback operators~\req{diagfdbk} respectively.
    \end{lemma}

These expressions can then be worked out for the variety of output operators $C$ representing the different performance measures defined earlier. The next result presents a summary of those calculations for the six different cases.

    \begin{corollary}
\label{var.thm}
The following are performance measures {\bf (P1)}, {\bf (P2)} and {\bf (P3)} expressed in terms of the Fourier symbols $\gh$, $\fh$ and $\ah$, of the operators $G$, $F$, and $T_a$ defining vehicular formations and consensus algorithms which satisfy assumptions {\bf (A1)}-{\bf (A5)}. The array $O$ is that of the standard consensus
algorithm~\req{scp_array}.
\begin{enumerate}
\item \emph{Consensus}
	\begin{enumerate}
		\item Local Error:
			\be
				V_{c}^{loc} ~=~\frac{1}{4d} \frac{1}{\beta} \sum_{n\neq 0,~ n\in\Z_N^d}
				\frac{ \Oh_n }{\Re(\ah_n)}
			   \label{Vcloc.eq}
			\ee		
		\item Long Range Deviation:
			\be
				V_{c}^{lrd} ~=~ -2
				\sum_{n_1 + \cdots + n_d~ \mbox{odd},~n\in\Z_N^d } \frac{1}{\Re(\ah_n)}
			   \label{Vclrd.eq}
			\ee
		\item Deviation from Average:
		\be
		V_c^{dav} ~=~ -\frac{1}{2} \sum_{n \neq0, ~ n\in\Z_N^d} \frac{1}{\Re(\ah_n)}.
		\label{cdav.eq}
		\ee
	\end{enumerate}

\item \emph{Vehicular Formations}
	\begin{enumerate}	
		\item Local Error:
			\be
			V_{v}^{loc} ~=~ -\frac{1}{4} \frac{1}{\beta}  \sum_{n\neq 0,~n\in\Z_N^d}
				\frac{ \Oh_n }{\gh_n \fh_n}
			   \label{vloc.eq}
			\ee
		\item Long Range Deviation:
		\be
			V_{v}^{lrd} =  2d \sum_{n_1 + \cdots + n_d ~\mbox{odd},~n\in\Z_N^d } \frac{1}{\gh_n \fh_n}
			\label{vlrd.eq}
		\ee
		\item Deviation from Average:
			\be
			V_{v}^{dav} = \frac{d}{2} \sum_{n \neq 0,~ n\in\Z_N^d} \frac{1}{\gh_n \fh_n}.
			\label{vdav.eq}
			\ee
	\end{enumerate}
\end{enumerate}
\end{corollary}

\section{Upper bounds using standard algorithms}
	\label{upbounds.sec}

In this section we derive asymptotic upper bounds for all three performance measures of
both the consensus and vehicular problems. These bounds are derived by
exhibiting simple feedback
laws similar to the one in the standard consensus algorithm~\req{scp_array}.
In the case of vehicular formations, we make a distinction between the cases of relative versus
absolute position and velocity feedbacks, and derive bounds for all four possible combinations of such
feedbacks.

The behavior of the asymptotic bounds has an important dependence on the underlying spatial
dimension $d$. For the purpose of cross comparison,
all of the upper bounds derived in this section are summarized in Table~\ref{upbounds.table}.

For later reference, we note that
the Fourier transform of the array $O$ in Eq.~\req{scp_array} is a quantity that occurs often, and can be easily calculated as
\begin{eqnarray}
	\Oh_n & = & -2d\beta ~+~ \sum_{r=1}^d  \left( \beta ~e^{-i \frac{2\pi}{N} n_r } ~+~  \beta ~e^{i \frac{2\pi}{N} n_r } \right)
					\nonumber  		\\
	  	& = & - 2\beta ~ \sum_{r=1}^d \left( 1 -  \cos\left(\frac{2\pi}{N}n_r   \right)  \right) .
   	             \label{ahsymb.eq}
\end{eqnarray}

\subsection{Upper bounds in the consensus case}

We consider the standard consensus algorithm~\req{cncsform}. In this case the array $a$ is
exactly $O$, and thus expression~\req{Vcloc} for the local error immediately simplifies to
\[
	V_c^{loc} ~=~   \frac{1}{4d\beta}  \sum_{n\neq 0,~ n\in\Z_N^d}  1
			~=~ \frac{1}{4d\beta} (M-1) ,
\]
which then implies the following upper bound for the individual local error at each site
\[
	\frac{ V_c^{loc} }{M} ~\leq~ \frac{1}{4d\beta}.
\]
Thus, the individual local error measure for the standard consensus
algorithm is bounded from above for any network size in any dimension $d$.

The derivation of the macroscopic error upper bounds are a little more involved.
First we observe that $V_c^{lrd}\leq 4 V_c^{dav}$. This is easily seen since
first, the sums in~\req{Vclrd} and~\req{cdav} involve terms that are all
of the same sign (since $\ah_n\leq 0$), and second, that the sum in~\req{Vclrd} is taken over
a subset of the terms in~\req{cdav}. It therefore suffices to derive the upper bounds for $V_c^{dav}$.

We begin with a
simplifying observation. Because the arrays $a$ we consider are real,
their Fourier symbols $\ah$  have even symmetry about all the mid
axes of $\Z^d_N$. More precisely
\[
	\ah_{(n_1,\ldots,n_r, \ldots,n_d)} ~=~ \ah_{(n_1,\ldots,N-n_r, \ldots,n_d)},
\]
for any of the dimension indices $r$. Assume for simplicity that $N$
is odd, and define $\Nb := (N+1)/2$. The even symmetry property implies that the
discrete hyper-cube $\Z_N^d$ can be divided into $2^d$ hyper-cubes, each of the size of $\Z_\Nb^d$, and
over which the values of $\ah$ can be generated from its values over $\Z_\Nb^d$ by appropriate reflections.
Consequently, a sum like~\req{cdav} can be reduced to
\[
		V_c^{dav}
		~ = ~ -\frac{1}{2} \sum_{n \neq0, ~ n\in\Z_N^d} \frac{1}{\Re(\ah_n)} ~=~
		-\frac{2^d}{2} \sum_{n \neq0, ~ n\in\Z_\Nb^d} \frac{1}{\Re(\ah_n)} .
\]

We now calculate an upper bound on the deviation from average measure~\req{cdav}
for the Fourier symbol~\req{ahsymb} of the
standard consensus algorithm
\begin{eqnarray}
	V_c^{dav} & = &  \frac{1}{4\beta}
		  \sum_{ n\neq 0, ~n\in\Z_N^d } \frac{1}{
			 \sum_{r=1}^d \left( 1 -  \cos\left(\frac{2\pi}{N}n_r   \right)  \right)	}	\nonumber	\\
	  & = &   \frac{2^d}{4\beta}  \sum_{ n\neq 0, ~n\in\Z_\Nb^d } \frac{1}{
			 \sum_{r=1}^d \left( 1 -  \cos\left(\frac{2\pi}{N}n_r   \right)  \right)	}	\nonumber	\\
	   & \leq & \frac{2^d}{32\beta} N^{2}
			\sum_{n\neq 0,~ n\in\Z_\Nb^d} \frac{1}{ (n^2_1+\cdots+n^2_d) }	,
	\label{cos2n2.eq}
\end{eqnarray}
where the first equality follows from reflection symmetry, and the inequality follows
from~\req{coslbound}, and noting that the
denominator is made up of $d$ terms of the form
\[
	1 ~-~ \cos\left( \frac{2\pi}{N} n_r\right) ~\geq~ \frac{2}{\pi^2} \left(\frac{2\pi}{N} n_r \right)^2
	~=~ \frac{8}{N^2} n_r^2,
\]
where the inequality is valid in the range $n_r \in [0,(\Nb-1)]$.

The asymptotics of sums in Eq.~\req{cos2n2}
are presented in Appendix~\ref{bounds.appen}. Using
those expressions, we calculate the individual deviation from average measure at each
site
\begin{eqnarray}
	\frac{V^{dav}_c}{N^d}
	& \leq &
		\frac{2^d}{32\beta} N^{2-d}
			\sum_{n\neq 0,~ n\in\Z_\Nb^d} \frac{1}{ (n^2_1+\cdots+n^2_d) }		 \nonumber \\
	& \approx &
			\frac{2^d}{32\beta} N^{2-d}
			\left\{ \begin{array}{ll}
			\frac{1}{d-2} (\Nb^{d-2}-1) ~~~~~ & d\neq 2	\\ \\
				\log(\Nb) 					& d = 2
				\end{array}		\right.			\nonumber	\\
	& \leq & 	
			{\cal C}_d \frac{1}{\beta}
			\left\{ \begin{array}{ll}
			N   & d=1 \\
				\log(N) 					& d = 2	\\
				1 ~~~~~ & d\geq 3
				\end{array}		\right.	,		
			\label{cncsasymp.eq}
\end{eqnarray}
where we have used $\Nb\leq N$, and ${\cal C}_d$ is a constant that depends on the dimension $d$,
but is independent of $N$ or the algorithm parameter $\beta$.
We note that the
upper bounds have exactly the same form when written in terms of
the network size $M=N^d$.

\begin{table*}[ttt]
	\caption{\textup{Summery of  asymptotic scalings of upper bounds in terms of the the total
		network size $M$ and the spatial dimensions $d$.
		Performance measures are classified as
		either microscopic (local error), or macroscopic (deviation from average or long range
		deviation). There are four possible feedback strategies in vehicular formations depending
		on which combination of relative or absolute position or velocity error feedback is used.
		Quantities listed are up to a multiplicative factor that is independent of $M$ or
		algorithm parameter $\beta$.
		}}
			\label{upbounds.table}
	\centering
	\begin{tabular}{|c||c|c|}
				\hline
	& \textbf{Microscopic} &  \textbf{Macroscopic} \\
				\hline
	\textbf{Consensus} &  ${1/\beta}$  &
		$ \frac{1}{\beta} \left\{ \begin{array}{ll}  M   & d=1 \\  \log(M)  & d = 2  \\ 1 ~~~~~ & d\geq 3
			  \end{array}		\right.	$
									 \\
				\hline
	\begin{tabular}	{rl}	\textbf{Vehicular Formations} &  \\
							Feedback type:	& 	abs. pos. \& abs. vel. \end{tabular}
			& $1/\beta$  & $1$  \\
				\hline
	\begin{tabular}	{rl}	\textbf{Vehicular Formations} &  \\
							Feedback type:	& 	rel. pos. \& abs. vel.\\
								{\em or}  & abs. pos. \& rel. vel. \end{tabular}
			 & $1/\beta$  &
		$ \frac{1}{\beta} \left\{ \begin{array}{ll}  M   & d=1 \\  \log(M)  & d = 2  \\ 1 ~~~~~ & d\geq 3
			  \end{array}		\right.	$			
			 						\\
				\hline
	\begin{tabular}	{rl}	\textbf{Vehicular Formations} &  \\
							Feedback type:	& 	rel. pos. \& rel. vel. \end{tabular}			&
				$ \frac{1}{\beta^2} \left\{ \begin{array}{ll}  M   & d=1 \\  \log(M)  & d = 2  \\ 1  & d\geq 3
			  		\end{array}		\right.	$				
			&
				$\frac{1}{\beta^2}
				\left\{ \begin{array}{ll}  M^3   & d=1 \\  M  & d = 2  \\ M^{1/3} & d=3 \\
									\log(M) & d=4	\\ 1 ~~~~~ & d\geq 5
			  		\end{array}		\right.	$							
			\\
	\hline	
	\end{tabular}
\end{table*}

\subsection{Upper bounds for vehicular formations}
	\label{upvf.sec}

To establish upper bounds in this case, we use a feedback control law which is similar
to~\req{platoonfdbk}. This law can be most compactly written in
operator notation as
\[
	u ~=~ T_O \xt ~+~T_O \vt ~+~ g_o \xt ~+~ f_o \vt,
\]
where $T_O$ is the operator of convolution with the array $O$ defined in the consensus
problem~\req{scp_array}. Note that in the multi-dimensional case, all signals are $d$-vectors,
and thus $T_O$ above is our notation for a diagonal operator with $T_O$ in each entry of the diagonal.
The last two terms represent absolute position and velocity error feedbacks respectively. The first two
terms represent a feedback where each vector component of $u_k$ is formed by a law
like~\req{cncsform} from
the corresponding vector components of $\xt_k$ and $\vt_k$ and all $2d$ immediate
neighbor sites in the lattice.

With the above feedback law, the closed loop system~\req{vfsys}
has the following expressions for the Fourier symbols of $G$ and $F$
\be
	\gh_n ~=~ g_o ~+~ \Oh_n , ~~~~~ \fh_n ~=~ f_o ~+~ \Oh_n,
   \label{gfubdef.eq}
\ee
where $\Oh$ is the Fourier symbol~\req{ahsymb}. We impose the additional conditions that
$g_o\leq 0$ and $f_o\leq 0$ since otherwise the closed loop system will have an increasing number
of strictly unstable modes as $N$ increases. When $g_o\neq 0$ (resp. $f_o\neq 0$) we refer
to that feedback as using absolute position (resp. velocity) feedback. There are four possible
combinations of such feedback scenarios.

We now use these expressions for the symbols $\gh$ and $\fh$ to calculate upper bounds on performance measures {\bf (P1)}, {\bf (P2)} and {\bf (P3)} for all four feedback scenarios. We begin with the local error~\req{vloc} which in this case is given by
    \be
	V_v^{loc} ~=~
	\frac{-1}{4\beta}    \sum_{n\neq 0,~n\in\Z_N^d}
				\frac{ \Oh_n }{(g_o+\Oh_n) (f_o+\Oh_n)}	 .
   \label{vlocbounds.eq}
    \ee
In the case of relative position and velocity error feedback, which corresponds to $g_o=0$ and $f_o=0$, the sum in Eq.~\req{vlocbounds} becomes $-\sum 1/\Oh_n$. This has the same form as $V_c^{dav}$ in Eq.~\req{cdav} for the standard consensus problem, and thus will grow asymptotically as derived in Eq.~\req{cncsasymp}. For this scenario, the final answer is listed as $V_v^{loc}$ in
Table~\ref{upbounds.table} after multiplying by the extra ${1/\beta}$ factor.
In the case of relative position and absolute velocity feedback,
the sum in Eq.~\req{vlocbounds} becomes $\sum -1/(f_o+\Oh_n)$. Each term is bounded from above by $-1/(f_o+\Oh_n) \leq -1/f_o$ since $f_o<0$ and $\Oh_n\leq 0$. Thus the entire sum has an upper bound that scales like $M$, which yields a constant bound for the individual local error once divided by the network size $M$.  An exactly symmetric argument applies to the case of absolute position but
relative velocity feedback. Finally, in the case of both absolute position and velocity feedback $f_o<0$ and $g_o<0$ implying a uniform bound on each term in the sum. Similarly the entire sum scales like $M$ and thus is uniformly bounded upon division by the network size. All of these four cases for the local error scalings are summarized in Table~\ref{upbounds.table}.

We now consider the case of the deviation from average measure~\req{vdav} which for our specific algorithm is
    \[
	V_{v}^{dav}
    \; = \; \frac{d}{2} \sum_{n \neq 0,~ n\in\Z_N^d} \frac{1}{(g_o+\Oh_n) (f_o+\Oh_n)}.
    \]
When $g_o<0$ and $f_o<0$, each term in the sum is bounded and the entire sum scales as $M$. Thus, the individual deviation from average at each site is bounded in this case. When either $f_o=0$ or $g_o=0$, then the
sums scale like $-\sum 1/\Oh_n$ (since the other factor in the fraction is uniformly bounded), i.e. like the deviation from average in the consensus case~\req{cncsasymp}.

The only case that requires further examination is that of relative position and relative velocity feedback
($g_o=f_o=0$). In this case
\begin{eqnarray*}
	V_{v}^{dav}
	& = &   \frac{d}{2} \sum_{n \neq 0,~ n\in\Z_N^d} \frac{1}{\Oh^2_n}       \\
	& \leq &   \frac{d2^d}{2^8} \frac{1}{\beta^2} N^{4}
			\sum_{n\neq 0,~ n\in\Z_\Nb^d} \frac{1}{ (n^2_1+\cdots+n^2_d)^2 },
\end{eqnarray*}
where the inequality is derived by the same argument used in deriving the inequality~\req{cos2n2}.
Dividing this expression by the network size $N^d$ and using the
asymptotic expressions~\req{multdsum2} yields
\be
	\frac{V^{dav}_v}{N^d}
	~ \leq ~
			{\cal C}_d \frac{1}{\beta^2}
			\left\{ \begin{array}{ll}
			\frac{1}{d-4} (1-N^{4-d}) ~~~~~ & d\neq 4	\\ \\
				\log(N) 					& d = 4
				\end{array}		\right. ,
   \label{upboundrelrel.eq}
\ee
where ${\cal C}_d$ is a constant depending on the dimension $d$ but independent of $N$ or the
algorithm parameter $\beta$.
Rewriting
these bounds in terms of the total network size $M=N^d$ gives the corresponding entries
in Table~\ref{upbounds.table}, where the other cases are also summarized.

We finally point out that $V_v^{lrd}\leq 4 V_v^{dav}$ due to an argument identical to that employed
in the consensus case. We thus conclude that the upper bounds just derived apply to the case of the
long range deviation measure as well.

\subsubsection*{The role of viscous friction}

It is interesting to observe that in
vehicular models with viscous friction~\req{wviscfric}, a certain amount of absolute velocity feedback
is inherently present in the dynamics. The model~\req{wviscfric} with a feedback control of the
form~\req{platoonfdbk} has the following Fourier symbol for the velocity feedback operator $F$
\[
	 \fh_n ~=~-\mu ~+~ f_o ~+~ \Oh_n.
\]
We conclude that even in cases of only relative velocity error feedback (i.e. when $f_o=0$),
the viscous friction term $\mu>0$ provides some amount of absolute velocity error feedback.
Thus, in an environment which has viscous damping, performance  in
vehicle formation problems scale  in a similar manner to consensus problems. These comments are
also applicable to the lower bounds developed in the next section.

\subsubsection*{The role of control effort}

A common feature of all the asymptotic upper bounds of the standard algorithms just presented is
their dependence on the parameter $\beta$.
If this parameter is fixed in advance based on
design considerations, then the algorithm's performance will scale as shown in
Table~\ref{upbounds.table}. However, it is possible to consider the redesign of the algorithms
as the network size increases. For example, it is possible to increase $\beta$ proportionally
to $M$ in consensus algorithms to achieve bounded macroscopic errors even for one dimensional
networks. As can be seen from~\req{cncsform}, this has the effect of increasing the
control feedback gains unboundedly (in $M$), which would clearly be unacceptable
in any realistic control problem. Thus, any consideration of the fundamental limits of
performance of more general algorithms must account for some notion of control effort,
and we turn to this issue in the next section.

\section{Lower bounds}
	\label{lbounds.sec}

A natural question arises as to whether one can design feedback controls with better asymptotic
performance than the standard algorithms presented in the previous section. In this section
we analyze the performance of any linear static state feedback control algorithm
satisfying the structural
assumptions {\bf (A1)}-{\bf (A5)}, and subject to a constraint on control effort.
A standard measure of control effort in  stochastic settings is the steady
state variance of the  control signal at each site
\be
	\expec{u^*_k u_k},
   \label{coneffdef.eq}
\ee
which is independent of $k$ due to the spatial invariance assumption. We constrain this
quantity and  derive lower bounds on the performance of any
algorithm that respects this constraint. The basic conclusion is that lower bounds on performance scale
like the upper bounds listed in Table~\ref{upbounds.table} with the control effort replacing the parameter
$\beta$.
In other words,  {\em any algorithm
with  control effort constraints will not do better
asymptotically than the standard algorithms} of section~\ref{upvf.sec}.
This is somewhat surprising given the extra degrees of freedom possible through feedback
control design, and it perhaps
implies that it is primarily the network topology and the structural constraints,
rather than the selection of the algorithm's parameters that determine these fundamental limitations.

We now turn to the calculation of lower bounds on both microscopic and macroscopic performance
measures. For brevity, we include only the calculations for the deviation from average macroscopic
measures.
These calculations are a little more involved than those for the upper bounds since
they need to be valid for an entire class of feedback gains. However, the basic ideas of
utilizing $H^2$ norms are similar, and this is what we do
in the sequel. In addition, a new ingredient appears where the control effort bound, combined with the
locality property, implies a uniform bound on the entries of the feedback arrays. This is stated
precisely in the next lemma whose proof is found in the Appendix.
These bounds then finally impose lower bounds on the performance of
control-constrained local algorithms.
\begin{lemma}
	\label{bounds.lemma}
Consider general consensus~\req{gen_cncs} and vehicular formation~\req{vfsys} algorithms
where the feedback arrays $a$, $g$ and $f$ posses the locality property (A3).
The following bounds hold
\begin{eqnarray}
	\|a\|_\infty & \leq & {\cal B}_a ~\expec{u^2_k} 	\nonumber    \\
	\|g\|_\infty & \leq & {\cal B}_g ~\left( \expec{u^2_k} \right)^2
								\label{arraybounds.eq}   \\
	\|f\|_\infty & \leq & {\cal B}_f ~\expec{u^2_k} 	,	\nonumber
\end{eqnarray}
where ${\cal B}_a$, ${\cal B}_g$ and ${\cal B}_f$ are constants independent of the network size.
\end{lemma}

\subsection{Lower bounds for consensus algorithms}
	
We start with the
deviation from average measure for a stable  consensus algorithm subject to a constraint
of  bounded control variance at each site
\be
	\expec{u_k^2} ~\leq~ W.
   \label{Wbound.eq}
\ee

We first observe a bound on $\Re(\ah_n)$ that can be established from the definition of the
Fourier transform
\begin{eqnarray*}	
	\Re(\ah_n)
	& \hspace*{-1em}= & 	  \hspace*{-1em}
		\Re\left(   \sum_{k\in\Z_N^d}  a_k e^{-i\frac{2\pi}{N}(n.k)} \right)
		=  \sum_{k\in\Z_N^d}  a_k \cos\left( \frac{2\pi}{N}n \cdot k \right)       \\
	& = & \sum_{k\in\Z_N^d}  a_k  \left[ 1 - \left( 1 -  \cos\left( \frac{2\pi}{N}n \cdot k \right)  \right)\right] \\
	& = & \sum_{k\in\Z_N^d}  (-a_k)  \left( 1 -  \cos\left( \frac{2\pi}{N}n \cdot k \right)  \right)
\end{eqnarray*}
where the last equality is a consequence of the condition $\sum_{k\in\Z_N^d} a_k = 0$. For lower bounds
on $\sum 1/\Re(-\ah_n)$, upper bounds on $\Re(-\ah_n)$ are needed. Observe that
\begin{eqnarray*}	
	| \Re(-\ah_n) |
	& = & \left| \sum_{k\in\Z_N^d}  a_k  \left( 1 -  \cos\left( \frac{2\pi}{N}n \cdot k \right)  \right)  \right|	\\
	& \leq &   \sum_{k\in\Z_N^d}  |a_k|  \left( 1 -  \cos\left( \frac{2\pi}{N}n \cdot k \right)  \right) 		 \\
	& \leq &	\frac{4\pi^2}{N^2} \sum_{k\in\Z_N^d}  |a_k|  ~ (n\cdot k)^2,
\end{eqnarray*}
where the second inequality follows from~\req{cosupbound}.
The last quantity can be further bounded by recalling the  locality property~\req{locality},
which has the consequence
\begin{eqnarray*}
	\lefteqn{
	 \sum_{k\in\Z_N^d} |a_k| ~ (k_1n_1+\cdots+k_dn_d)^2	} ~~~~~~~~~~~~~~~~	\\
	 & = &
	 \sum_{k\in\Z_N^d, ~|k_i|\leq q} |a_k| ~ (k_1n_1+\cdots+k_dn_d)^2	\\
	 & \leq &
	    \sum_{0\neq k\in\Z_N^d, ~|k_i|\leq q} |a_k| ~ (qn_1+\cdots+qn_d)^2	\\
	    & = &  q^2 (n_1+\cdots+n_d)^2  \sum_{0\neq k\in\Z_N^d, ~|k_i|\leq q} |a_k|  .
\end{eqnarray*}
Now the locality property can be used again to bound the above sum using the
the control effort bounds~\req{arraybounds} and~\req{Wbound}
\be
	\sum_{0\neq k\in\Z_N^d, ~|k_i|\leq q} |a_k|
	 ~ \leq ~
		(2q)^d \|a\|_\infty ~\leq~ (2q)^d  ~\cB_a ~W  .
   \label{l1Wbd.eq}
\ee
Putting the above together gives
\begin{eqnarray*}
	V_c^{dav}
	& = &  \sum_{n\neq 0, ~n\in\Z_N^d} \frac{1}{-\Re(\ah_n)} 	\\	
	& \geq &  \frac{N^2}{\pi^2  (2q)^{d+2} ~\cB_a W}
			\sum_{n\neq 0,~ n\in\Z_N^d}
			~\frac{1}{ (n_1+\cdots+n_d)^2} 			\\
		& \geq &   \frac{{\cal C}_d}{W} ~N^2 \sum_{n\neq 0,~ n\in\Z_N^d}
			~\frac{1}{ (n^2_1+\cdots+n^2_d)} ,
\end{eqnarray*}
where the last inequality follows from~\req{n2bound}, and ${\cal C}_d$ is a constant
independent of $N$.

Finally, utilizing~\req{multdsum} and dividing by the network size $M=N^d$,  a lower
bound on the deviation from average is obtained
\bea
	\frac{V^{dav}_c}{N^d}
	& \geq & \frac{{\cal C}_d}{W} ~N^{2-d} \sum_{n\neq 0,~ n\in\Z_N^d}
			~\frac{1}{ (n^2_1+\cdots+n^2_d)}				\nonumber	\\
	& \approx &
			 \frac{{\cal C}_d}{W}
			 \left\{ \begin{array}{ll}
				\frac{1}{d-2} (1-N^{2-d}) ~~~~~ & d\neq 2		\\ \\
				\log(N) 					& d = 2
				\end{array}		\right. ,				\nonumber	\\
	& \geq & 	
			{\cal C}_d \frac{1}{W}
			\left\{ \begin{array}{ll}
			N   & d=1 \\
				\log(N) 					& d = 2	\\
				1 ~~~~~ & d\geq 3
				\end{array}		\right.	,
	\label{cons_result.eq}
\eea
where by a slight abuse of notation, we use ${\cal C}_d$ to denote different constants in the
expressions above.  We observe how the lower bounds~\req{cons_result} have the same
asymptotic form as the upper bounds for the standard consensus algorithm~\req{cncsasymp},
but with the control effort bound $W$ replacing the parameter $\beta$.

\subsection{Lower bounds for vehicular formations}

We recall the development of the upper bounds for vehicular formations in Section~\ref{upvf.sec}.
The Fourier symbols of general feedback gains $G$ and $F$ have a similar form to~\req{gfubdef},
and can be written as
\be
	\gh_n ~=~ g_o ~+~ \gammah_n , ~~~~~ \fh_n ~=~ f_o ~+~ \phih_n,
   \label{gflbdef.eq}
\ee
where $g_o$, $f_o$ and $\gammah$, $\phih$ are the absolute and relative feedback terms respectively.
As before, we impose the conditions that $g_o,f_o\leq 0$. We assume that we have a control effort
constraint of the form~\req{Wbound}.

The case of absolute position and absolute velocity feedback has upper bounds which are finite, and
the question of lower bounds is moot. For the other three cases, lower bounds on~\req{vdav}
are established using upper bounds on the symbols $\gh$ and $\fh$ which can be derived as follows
\[
	\|\fh\|_\infty ~\leq~ \|f\|_1 ~\leq~ (2q+1) ~\|f\|_\infty ~\leq~ (2q+1)\cB_f ~W ,
\]
where the inequalities follow from~\req{inftybound}, the locality property,
and~\req{arraybounds} respectively. For $g$ we similarly have
\[
	\|\gh\|_\infty ~\leq~ (2q+1)\cB_g ~W^2 .
\]

Consider now the case of relative position and absolute velocity feedback. A lower bound is
established by
\begin{eqnarray*}
	V_{v}^{dav} &=&  \frac{d}{2} \sum_{n \neq 0,~ n\in\Z_N^d} \frac{1}{|\gh_n| |\fh_n|}		\\
			&\geq&  \frac{d}{2(2q+1)\cB_f}	~\frac{1}{W}
					\sum_{n \neq 0,~ n\in\Z_N^d} \frac{1}{|\gh_n|}  .
\end{eqnarray*}
Now a lower bound on the sum can be established in exactly the same manner as~\req{cons_result}
in the consensus case since $\gh$ is a symbol of a local relative feedback operator. The case of
relative velocity and absolute position feedback is similar with the exception that the factor of
$\frac{1}{W}$  is replaced by $\frac{1}{W^2}$.

The final case to consider is that of relative position and relative velocity feedback. One can
repeat the same arguments made in the consensus case up to equation~\req{l1Wbd} for both
$\gh_n$ and $\fh_n$ to state
\begin{eqnarray*}
 	\sum_{0 \neq n\in\Z_N^d} \frac{1}{|\gh_n| |\fh_n|}
	&\geq &      \frac{c_1 ~~N^4}{\|g\|_\infty \|f\|_\infty}
		 \sum_{0 \neq n\in\Z_N^d}  \frac{1}{ (n_1+\cdots+n_d)^4}			\\
	&\geq &   \frac{c_2 ~~N^4}{W^3}
		 \sum_{0 \neq n\in\Z_N^d}  \frac{1}{ (n_1+\cdots+n_d)^4}	,
\end{eqnarray*}
where $c_1$ and $c_2$ are some constants independent of $N$ and $W$. The asymptotic
behavior of this expression (divided by the network size)
was given earlier in~\req{upboundrelrel}. We thus conclude that the
lower bounds in this case are exactly like the upper bounds shown in Table~\ref{upbounds.table}
for relative position and relative velocity feedback,
but with the $\frac{1}{\beta^2}$
term replaced by $\frac{1}{W^3}$.

\section{Examples and Multiscale Interpretation}
	\label{comp.sec}

\begin{figure}[h]
	\begin{center} 	\includegraphics[width=3in]{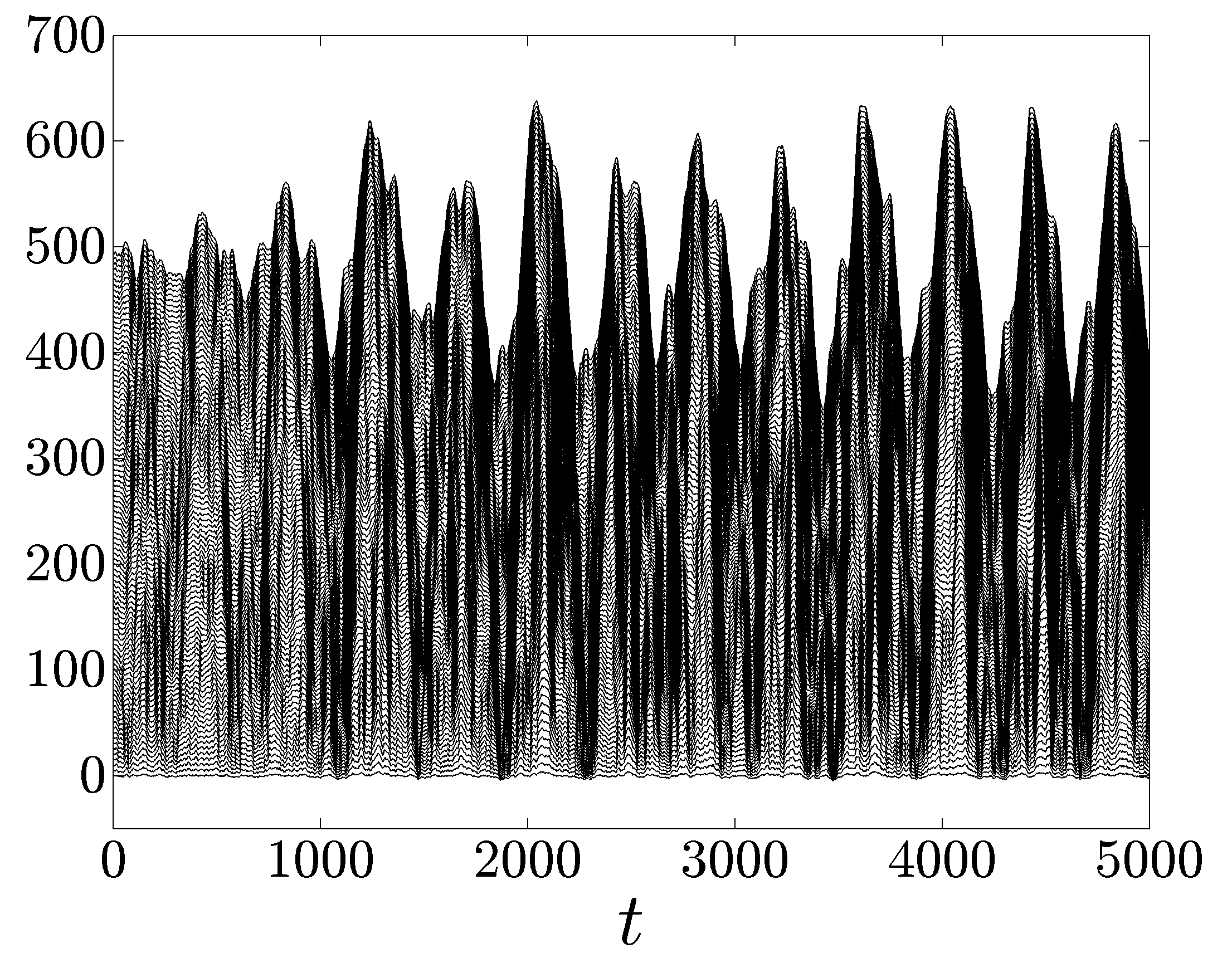}  \end{center}
	\begin{center} 	\includegraphics[width=2in]{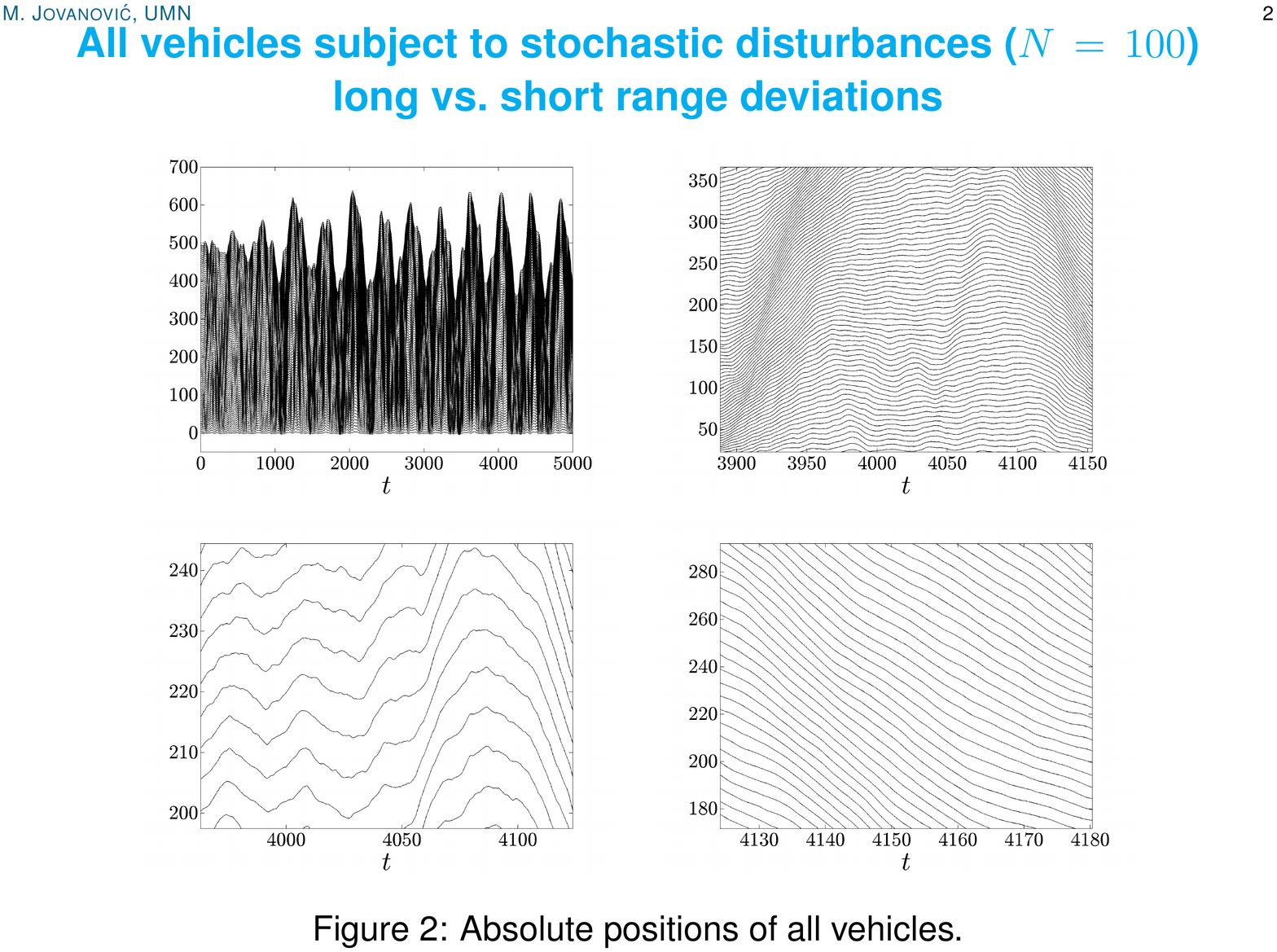}  \end{center}
	\caption{Vehicle position trajectories (relative to vehicle number 1) of a 100 vehicle formation all of which
		are subjected to random disturbances. Top graph is a  ``zoomed out'' view
		exhibiting the slow accordion-like motion
		of the entire formation. Bottom graph is a zoomed in view showing that vehicle-to-vehicle
		distances are relatively well regulated.}
   \label{accordion.fig}
\end{figure}

Numerical simulations of cases
where macroscopic measures grow unboundedly with network size show a particular type of motion
for the entire formation. In the one dimensional case, it can be described as an accordion-like motion
in which large shape features in the formation fluctuate. Figure~\ref{accordion.fig} shows the results
of a simulation of a 100 vehicle platoon with both relative position and relative velocity error
feedbacks. This corresponds to a control strategy of the type for which upper bounds were
calculated in section~\ref{upvf.sec} with $g_o=f_o=0$.

An interesting feature of these plots is the phenomenon of lack of formation coherence.
This is only
discernible when one ``zooms out'' to view the entire formation.
The length of the formation
fluctuates stochastically, but with a distinct slow temporal and long spatial wavelength signature.
In contrast,
the zoomed-in view in Figure~\ref{accordion.fig} shows a relatively well regulated vehicle-to-vehicle
spacing. In general, it appears that small scale (both temporally and spatially) disturbances are
well regulated, while large scale disturbances are not. An intuitive interpretation of this phenomenon
is that local feedback strategies are unable to regulate against large scale disturbances.

In this paper we have not directly analyzed the temporal and spatial scale dependent disturbance
attenuation limits of performance. However, it appears that our microscopic and macroscopic
measures of performance do indeed correspond to small and large scale (both spatially and
temporally) motions respectively. We next outline a more mathematical argument that connects
these measures.

\subsubsection*{Mode shapes}
To appreciate the connection between $H^2$ norms and mode shapes in our system, consider
first a general  linear system driven by a white  random process
\[
	\dot{x} ~=~ Ax ~+~ w.
\]
When $A$ is a normal matrix, it is easy to show (by diagonalizing the system with the orthonormal
state transformation made up of the eigenvectors of $A$) that the steady state variance of the
state is
\[
	\lim_{t \rightarrow \infty}	{\cal E}\{x^*(t)x(t)\} ~=~ \sum_i \frac{1}{2~\Re\{\lambda_i\}},
\]
where the sum is taken over all the eigenvalues $\lambda_i$ of $A$. Thus we can say that under
white disturbance excitation, the amount of energy each mode contains is inversely proportional to its
distance from the imaginary axis. In other words, slower modes are more energetic. Now, all the systems
we consider in this paper are diagonalizable (or block-diagonlizable) by the spatial Discrete Fourier
Transform. In addition, for the standard algorithms, we have the situation that
{\em slow temporal modes correspond to long spatial
wavelengths}.
This provides an explanation for the observation that the most energetic motions are those
that are temporally slow and have
long spatial wavelengths.

\subsubsection*{String instability}
\begin{figure}[h]
	\begin{center} 	\includegraphics[width=2in]{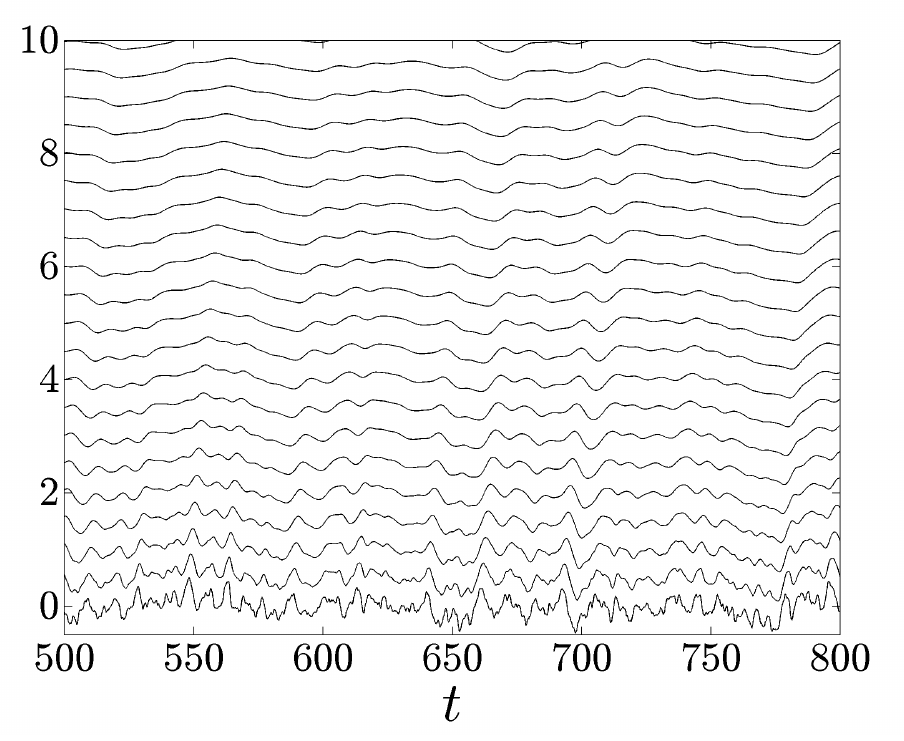}  \end{center}
	\caption{Vehicle position trajectories (relative to leader) of the first few of a 100 vehicle
		formation. Only lead vehicle is subjected to random disturbances. Vehicle trajectories
		exhibit regulation against that disturbance, indicating the absence of string instability. }
   \label{nsi.fig}
\end{figure}

While string instability is sometimes an issue in formation control, the phenomenon we study
in this paper is distinct from string instability. The example presented in this section
 is that of a formation
that does posses string stability. For illustration, we repeat the simulation but with disturbances
acting only on the first vehicle. The resulting vehicle trajectories are shown in Figure~\ref{nsi.fig}.
It is interesting to note that temporally high frequency disturbances appear to be very well regulated,
and do not propagate far into the formation, while low temporal frequency disturbances
appear to propagate deep into the formation. What is not shown in the figure is that low frequency
disturbances are eventually  regulated for vehicles far from the first. This is consistent with
the intuitive notion discussed earlier that local feedback is relatively unable to regulate large scale
disturbances.

\subsubsection*{Multi-scale properties of disturbance rejection}
An intriguing explanation of the above example and our scaling results is as follows. The macroscopic
error measures capture how well the network regulates against large-scale disturbances. In large, one
dimensional networks, local feedback alone is thus unable to regulate against these large-scale
disturbances, and global feedback is required to achieve this. This seems rather intuitive.
Perhaps surprisingly, in large networks with
higher spatial dimensionality, local feedback alone can indeed regulate against large-scale disturbances.
This follows for networks for which the macroscopic error measure is bounded irrespective of network
size.
The ``critical dimension'' needed to achieve this depends on the order of the node dynamics as well as the
type of feedback strategy as shown in Table~\ref{upbounds.table} (e.g. dimension 3 for relative position
and absolute velocity feedback, and dimension 5 for relative position and velocity feedback in cases of
vehicular formations).

\section{Discussion}
	\label{disc.sec}

\subsection{General networks}

The networks considered in this paper are ones which can be built on top of a Torus network. Some
concepts, such as coherence and microscopic and macroscopic errors are easily generalized to
arbitrary networks. The correct generalization of the concept of spatial dimension however is more
subtle.

For any network of dynamical systems for which a distance metric is defined between nodes (e.g. from
an imbedding of the network in $\R^n$), the notion of long range deviation can be defined as done in this
paper. The calculation of that quantity involves system Grammians and may even be written in terms of
the underlying system matrices for certain structures. Thus coherence  measures
can be calculated numerically for
such networks. However, more explicit calculations to uncover scaling laws as network size increases
will clearly require more analytical expressions for the system norms in such networks.

To generalize the present results, one would require a notion of how to grow the network size while
preserving certain topological properties such as the spatial dimension.
 Preliminary results on
self-similar and fractal networks have been obtained~\cite{patbam11cdc}.
The proper notion of spatial dimension to capture coherence in  general
graphs  remains a research topic at this time.

\subsection{Distributed estimation and resistive lattices}

The results presented here have a strong resemblance to results on performance limitations
of distributed estimation algorithms based on network topology~\cite{barhes07,barhes08},
where asymptotic bounds similar to~\req{cncsasymp} first appeared in the controls literature
(see also~\cite{barhes06cdc} where a consensus problem with noisy observations is analyzed yielding
performance bounds like the consensus upper bounds we have in the present paper).
In that work, the arguments are based
on an analogy with effective resistance in resistive lattices and certain imbeddings of their graphs
in d-dimensional space~\cite{doyle1984rwa}.
It is not clear how the resistive analogy can be generalized to cover the case of
second order dynamics (i.e. vehicular formations), or the lower bounds on more general control laws.
We have therefore avoided the resistive network analogy in this paper
by directly using the multi-dimensional Fourier transform and reducing all calculations
to sums of the form~\req{multdsum} resulting in a self-contained argument.

It is interesting to note that the original arguments
for the asymptotic behavior in resistive lattices~\cite{cserti2000alg} in the physics literature are based
on approximations of the Green's function of the diffusion operator in $d$-dimensions, for which
the underlying techniques are approximations like~\req{integapp}.

\subsection{Order of local dynamics}

We have attempted to keep the development general enough that it is applicable to networked
dynamical systems whose dynamics are not necessarily those of vehicles in any particular
physical setting.
What we refer to in this paper as consensus and vehicular formations problems respectively
represent networks where the local dynamics (at each site) are first and second order chains of
integrators respectively.
The dynamical models are such that the stochastic disturbance enters into the first integrator, and
the performance objectives involve variances of the outputs of the last integrators at each site.
One generalization of this set up is where the local dynamics is a chain of $n$ integrators. It is then
possible to show that (by retracing the arguments for the vehicular formations case and
generalizing~\req{multdsum2})
the cutoff dimension to have bounded macroscopic measures with only local relative state feedback
is $1+2n$.

\subsection{LQR designs}

It was observed in~\cite{jovbamTAC05platoons} that optimal LQR designs for vehicular platoons suffer
from a fundamental problem as the platoon size increases to infinity. These optimal feedback
laws are almost local in a sense described by~\cite{bampagdah02}, where control gains
decay exponentially as a function of distance. The resulting optimal feedbacks~\cite{jovbamTAC05platoons}
suffer from the
problem of having underdamped slow modes with long spatial wavelengths. Thus, the same incoherence
phenomenon occurs in these optimal LQR designs where the performance objective is composed of
sums of local relative errors (leading to  feedback laws with exponentially decaying gains on relative
errors).

\subsection{Measuring performance in large scale systems}

	In spatial dimensions where performance scalings are bounded,
	the underlying system eigenvalues still
	limit towards zero, suggesting ultimate instability in the limit as $M \rightarrow \infty$. However,
	measures of performance remain bounded in these cases. In such cases the locations
	of internal eigenvalues are not a good indication of the system's performance in the limit of
	large networks.
	
	Take the consensus problem over $\Z_N^d$ as an example. The ``least damped eigenvalue''
	(other than zero)
	quantifies the convergence time of deviation from average (in the absence of
	stochastic disturbances), and it scales as
	\be
		\frac{1}{|\lambda_2|} ~=~ \Theta\left({N^{2/d}}\right) ,
	   \label{lde.eq}
	\ee
	as can be shown by explicit eigenvalue calculations~\cite{carli2007cca, PBE06}.
	If we use this quantity as a measure of performance, it indicates that performance
	becomes arbitrarily bad (in the limit of large $N$) in any spatial dimension $d$.
	On the
	other hand, consider the use of a macroscopic error measure like the variance of the deviation
	from average~\req{davcons} in the presence of stochastic disturbances.
	That quantity can be expressed in terms of the system
	eigenvalues as
	\be
		V_c^{dav} ~=~ \frac{1}{2N^d}  \sum_{n\neq 0} \frac{1}{|\lambda_n|},
	   \label{davassum.eq}
	\ee
	where the sum is taken over all the system's eigenvalues other than zero. Note
	that this sum is just~\req{cdav} rewritten to emphasize the contrast with~\req{lde}.
	
	The important observation is that~\req{lde} indicates that as network size increases, the
	system eigenvalues approach the stability boundary, indicating an eventual catastrophic
	loss of performance in any spatial dimension $d$. On the other hand, ~\req{davassum} is
	uniformly bounded in dimensions $d\geq 3$ (as shown in~\req{cncsasymp}),
	thus implying well behaved
	systems as quantified by the macroscopic performance measures. A similar point to the
	above has been recently made~\cite{cargarzam09ita}.
	
	The least damped eigenvalue is traditionally used as an important measure of performance.
	The examples in this paper demonstrate that for large scale systems, it is not a very
	meaningful measure of performance, and that the general question of how to measure
	performance in large scale systems is a subtle one.

\subsection{Detuning/mistuning designs}

It is shown in~\cite{barmehhes09} that spatially-invariant  local controllers for platoons
have closed loop eigenvalues that approach the origin at a rate of $O(\frac{1}{N^2})$. A
``mistuning'' design modification is proposed~\cite{barmehhes09}, resulting in spatially-varying
local controllers where the closed loop eigenvalues approach the origin at the better rate of
$O(\frac{1}{N})$. In this paper, we have not used the real part of the least damped eigenvalue
as a measure of performance but rather the variance of certain system outputs. This amounts
to using an $H^2$ norm as the measure of performance. It was shown in ~\cite{bampagdah02}
that for spatially-invariant plants, one can not improve $H^2$ performance with spatially-varying
controllers. The resulting controllers however have exponentially decaying gains rather than
completely local gains. The problem of designing optimal $H^2$ controllers with a prescribed
neighborhood of interaction remains an open and non-convex one.
It is an interesting and open question as to whether mistuning designs for the $H^2$ measures
we use in this paper
can yield local controllers with better asymptotic performance than spatially-invariant ones.
It was also shown~\cite{barmehhes09} that a mistuning design can improve $H^\infty$ performance
for platoon problems. This shows that there is perhaps an important distinction between $H^\infty$ and
$H^2$ measures of performance for large scale systems. A point that is worthy of further investigation.

\appendix

\subsection{Multi-dimensional Discrete Fourier Transform}
	\label{fourier.appen}

We define the Discrete Fourier Transform for functions $\{f_k\}$ over $\Z_N^d$ by
\[
	\fh_n ~:=~ \sum_{k\in\Z_N^d} f_k ~e^{-i\left(\frac{2\pi}{N} n\cdot k \right) },
\]
where $n\cdot k ~:=~ n_1k_1+\cdots + n_dk_d$. The inverse transform is given by
\[
	f_k ~:=~  \frac{1}{M} \sum_{n\in\Z_N^d} \fh_k ~e^{i\left(\frac{2\pi}{N} n\cdot k \right) },
\]
where $M=N^d$. An immediate consequence of the definitions are the following bounds
\be
	\|\fh\|_\infty ~\leq~ \|f\|_1 ,
	~~~
	\|f\|_\infty ~\leq~ \frac{1}{M} \|\fh\|_1 .
   \label{inftybound.eq}
\ee

Let $\delta$ be the Kronecker delta on $\Z_N^d$. It's transform is the array $\bfo$, which is the
array of all elements equal to $1$. The transform of $\bfo$ is $M\delta$.
We use the symbol $\bfob$ to denote 	the array of all elements equal to $\frac{1}{M}$.

If $T_f$ denotes the circulant operator of circular convolution with $f$, then the eigenvalues of
$T_f$ are just the numbers $\{\fh_n\}$, and consequently the trace of $T_f$ is given by the sum
\be
	\mbox{tr}\left( T_f \right) ~=~ \sum_{n\in\Z_N^d} \fh_n.
   \label{circtr.eq}
\ee

\subsection{Bounds and asymptotics of sums}
	\label{bounds.appen}

		The following facts are useful in establishing asymptotic bounds.
		
		\noindent	\textbf{1)} For any $x\in\R$ and any $y\in[-\pi,\pi]$
			\begin{eqnarray}	
				1-\cos(x) & \leq & x^2,
			 	\label{cosupbound.eq}				\\					
				1-\cos(y)  & \geq &  \frac{2}{\pi^2}~ y^2.
			 	\label{coslbound.eq}	
			\end{eqnarray}
		\noindent	\textbf{2)}
			Given $d$ integers $n_1$, $\ldots$, $n_d$,
			\be
				(n_1+\cdots+n_d)^2 ~\leq~ (2d+1)~(n_1^2+\cdots+n_d^2).
			   \label{n2bound.eq}
			\ee
				\textbf{Proof:}{
				\[
							\left(\sum_i n_i \right)^2 ~ = ~ \sum_i n_i^2 ~+~
										\sum_{i} \sum_{j\neq i} n_i n_j .
				\]
				Using $n_i n_j ~\leq~ \left( \max\{n_i,n_j\}\right)^2 ~\leq~ n_i^2+n_j^2$, we get
				the bound
				\beas
					 \left(\sum_i n_i \right)^2
					 & \leq &  \sum_i n_i^2 ~+~
										\sum_{i=1}^d \sum_{j\neq i} (n_i^2 +n_j^2)	 \\
					& \leq &  \sum_i n_i^2 ~+~ 2d \sum_i n_i^2.
				\eeas
				}

		\noindent	\textbf{3)}  In the limit of large $N$,
			\begin{eqnarray}
			 \sum_{ \stackrel{\scriptstyle n\neq 0}{ n\in\Z_N^d} }  ~\frac{1}{ (n^2_1+\cdots+n^2_d)}
			  & \approx &
				\left\{ \begin{array}{ll}
					\frac{1}{d-2} (N^{d-2}-1) ~ & d\neq 2	\\ \\
					\log(N) 					& d = 2
					\end{array}		\right.
				   \label{multdsum.eq}						\\
			 \sum_{ \stackrel{\scriptstyle n\neq 0}{ n\in\Z_N^d} }  ~\frac{1}{ (n^2_1+\cdots+n^2_d)^2}
			  & \approx &
				\left\{ \begin{array}{ll}
					\frac{1}{d-4} (N^{d-4}-1) ~ & d\neq 4	\\ \\
					\log(N) 					& d = 4
					\end{array}		\right. 		
				   \label{multdsum2.eq}						\\							
			\end{eqnarray}
			where $f(N)  ~\approx~ g(N)$ is notation for
			\[
				\underbar{c} ~g(N) ~\leq~  f(N)  ~\leq~ \bar{c} ~g(N),
			\]
			for some constants $\bar{c}$ and $\underbar{c}$ and all $N\geq \bar{N}$ for some $\bar{N}$.
							\\
				\textbf{Proof:} {
				We begin with~\req{multdsum}.
				Upper and lower bounds on this sum can be derived by
				viewing it as upper and lower Rieman sums for the integral
				\[
					\int \cdots \int \frac{1}{x_1^2+\cdots+x_d^2} dx_1 \cdots dx_d,
				\]
				over the region $\Delta\leq r\leq 1$ for the lower bound, and
				$\Delta \leq r \leq \sqrt{d}$ for the upper bound. Here $\Delta=\frac{1}{N}$,
				and the asymptotic behavior is determined by the lower limit on the integrals,
				so both upper and lower bounds behave the same asymptotically.
				
				Using the grid points $x_1 = n_1\Delta$, $\ldots$, $x_d=n_d\Delta$, and
				using the volume increment $\Delta^d$, we get
				\begin{eqnarray}
					\lefteqn{
					\int \cdots \int \frac{1}{x_1^2+\cdots+x_d^2} dx_1 \cdots dx_d} \nonumber  \\
					& \approx &
					\Delta^d  \sum_{n\neq 0,~ n\in\Z_N^d}
						~\frac{1}{ ((\Delta n_1)^2+\cdots+(\Delta n_d)^2)}							 \nonumber	\\
					& = &
					\Delta^{d-2} \sum_{n\neq 0,~ n\in\Z_N^d}   ~\frac{1}{ (n^2_1+\cdots+n^2_d)}.
					   \label{integapp.eq}
				\end{eqnarray}
				Now the integral can be evaluated using hyperspherical coordinates by
				\begin{eqnarray*}
				\lefteqn{
				\int  \frac{1}{x_1^2+\cdots+x_d^2} dx_1 \cdots dx_d  ~=}     \\
				& &  \hspace*{-2em}  \int_\Delta^1 \int_0^\pi \cdots \int_0^{2\pi}
					\frac{1}{r^2} ~r^{d-1} \sin^{d-2}(\phi_1)\ldots\sin(\phi_{d-2}) ~dr~d\phi_1\ldots
								d\phi_{d-1}								 	\\
				 & & ~= C_d \int_\Delta^1 r^{d-3} dr,
				\end{eqnarray*}
				where $C_d$ is a constant that depends only on the dimension $d$ (and can be
				expressed in terms of the volume of the unit sphere in $\R^d$). Evaluating this
				integral, using $\Delta=\frac{1}{N}$ and equation~\req{integapp} gives the
				result~\req{multdsum}.
				
				The proof of~\req{multdsum2} is very similar to the above, with the exception
				that one approximates the integral of $\frac{1}{(x_1^2+\cdots+x_d^2)^2} =
				\frac{1}{r^4}$ instead. The details are omitted for brevity.
				}

	\subsection{Proof of Lemma \ref{V.lemma}}

		For the consensus problem, the state equation is~\req{gen_cncs}, and thus the Lyapunov
		equation~\req{pklyap} becomes
		\[
			\ah_n^* \ph_n ~+~\ph_n \ah_n ~=~ -\Ch_n^* \Ch_n,
		\]
		where we have used $H=C$ (and the choice of $C$ depends on the particular performance
		measure being considered). Since all quantities are scalars, this equation is immediately solved
		for $\ph_n ~=~ \frac{-\Ch_n^* \Ch_n}{2\Re(\ah_n)}$ for $n\neq 0$. In the case $n=0$, we look
		at the integral definition~\req{pkinteg}, conclude that $\Ch_0=0$ implies that $\ph_0=0$. Thus,
		the sum~\req{totalh2} is calculated to be~\req{Vc}.

		For the vehicular problem, the state equation is~\req{vfsys} with
		the output equation $H=\obt{C}{0}$.
		The Lyapunov equation~\req{pklyap} becomes
		\beas
			\lefteqn{
			\tbt{0}{\Gh_n^*}{I}{\Fh_n^*} \tbt{\Xh_n}{\Zh_n}{\Zh_n^*}{\Yh_n} +
			\tbt{\Xh_n}{\Zh_n}{\Zh_n^*}{\Yh_n} \tbt{0}{I}{\Gh_n}{\Fh_n} 			 }
			~~~~~~~~~~~~~~~~~~~~~~~~~~~~~~~~~
			&  & 			\\
			& =  &  \tbt{-\Ch_n^*\Ch_n}{0}{0}{0},
		\eeas
		where each of the submatrices is of size $d\times d$.
		From the above, we extract the following matrix equations
		\begin{eqnarray}
			\Gh_n^* \Zh_n^* ~+~\Zh_n \Gh_n  & = & -\Ch_n^*\Ch_n		\label{1of3.eq}	 \\
			\Gh_n^*\Yh_n ~+~\Xh_n ~+~\Zh_n \Fh_n & = & 0 		\nonumber	\\
			\Zh_n ~+~\Fh_n^*\Yh_n ~+~\Zh_n^* ~+~\Yh_n \Fh_n &=& 0 .	 \label{3of3.eq}
		\end{eqnarray}
		Since we are only interested in the quantity
		\[
			\mbox{tr}\left( \Bh_n^*\Ph_n \Bh_n \right) ~=~
			\obt{0}{I} \tbt{\Xh_n}{\Zh_n}{\Zh_n^*}{\Yh_n} \tbo{0}{I} ~=~
			\mbox{tr}\left( \Yh_n \right) ,
		\]
		then only equations~\req{1of3} and~\req{3of3} are relevant. The coordinate decoupling
		assumption (A5) on the operators
		$G$, $F$ and $C$ implies that the matrices $\Gh_n$, $\Fh_n$ and $\Ch_n$ are all diagonal.	
		It follows that $\Zh_n$, $\Xh_n$ and $\Yh_n$ are also
		diagonal, and the above matrix equations are trivial to solve.  $\Zh_n$ has the
		solution
		\[
			Z~=~ -\frac{1}{2} \Gh_n^{-1} \Ch_n^*\Ch_n.
		\]
		Similarly, equation~\req{3of3} is solved
		to yield
		\[
			Y ~=~ \frac{1}{2} (\Gh_n \Fh_n)^{-1} \Ch_n^*\Ch_n  ,
		\]
		for $n\neq 0$.
		For the unstable
		mode at $n=0$, the integral~\req{pkinteg} can be easily evaluated to yield $\Zh_0=0$ (since
		$\Ch_0=0$ for all the performance measures we consider).
		Adding in the assumption that all matrices are diagonal with equal elements, we obtain
		in summary the total $H^2$ norm
		of the vehicle formation problem~\req{vfsys} is given by
		\be
			 V_v ~=~ 	 \frac{d}{2} \sum_{n\neq 0~ n\in\Z_N^d} \frac{\ch^*_{n}\ch_{n}}{\gh_n\fh_n},
		\ee
		where the multiplicative factor of $d$ comes from taking the trace of $\Yh_n$.

	\subsection*{Proof of Corollary \ref{var.thm}}

		\subsubsection{Consensus}
	
		The local error measure output operator $C$ is given by~\req{Clocaldef}, for which
		$C^*C ~=~\frac{-1}{2d\beta} ~O$ by the identity~\req{Oiden}.
		Combining this with Lemma~\ref{V.lemma} gives the result for $V_c^{loc}$.
		
		The long range deviation measure has the output operator defined in~\req{Clrd}, which has
		the Fourier symbol
		\[
			\ch_n ~=~ 1 ~-~ e^{-i\pi \left(n_1+\cdots+n_d\right)},
		\]
		from which we conclude that
		\[
			\ch_n ~=~  \left\{ \begin{array}{lcl} 0 & ~~& \left(n_1+\cdots+n_d\right) ~\mbox{even} \\
												                 2 & & \left(n_1+\cdots+n_d\right) ~\mbox{odd}
										 \end{array} \right.   .
		\]			
		Combining this with Lemma~\ref{V.lemma} gives the result for $V_c^{lrd}$.

		The deviation from average output operator is $C ~=~ I -T_\bfo$, or equivalently,
		the convolution operator $T_{\deltab^0-\bfo}$. The corresponding Fourier symbol
		is the Fourier transform of the array $\deltab^0-\bfo$, which gives
		\[
			\ch_n ~=~ \bfo_n - \delta_n ~=~ \left\{ \begin{array}{lcl} 0 & ~~& n=0 \\
												                 1 & & n\neq 0
										 \end{array} \right.   .
		\]
		Putting this in the general formula~\req{Vc} yields the result for $V_c^{dav}$

	\subsubsection{Vehicular formations}
		The derivations for this case are very similar to those for the consensus problem and are therefore
		omitted for brevity.

	\subsection*{Proof of Lemma \ref{bounds.lemma}}

		We rewrite the dynamics of the consensus algorithm so that $u$ is an output
		\begin{eqnarray*}
			\dot{x} & = &  T_a x ~+~ w \\
			u & = &  T_a x.
		\end{eqnarray*}
		The present task is then to calculate
		the $H^2$ norm from $w$ to $u$. 		
		Applying formula~\req{Vc} with $T_a$ as the $C$ operator yields
		\[
			\sum_{k\in\Z_N^d} \expec{u^2_k}
				= \frac{-1}{2} \sum_{n\neq 0, ~n\in\Z_N^d} \frac{\ah_n^2}{\ah_n}
			    	~=~ \frac{1}{2} \sum_{n\in\Z_N^d}  (-\ah_n) ,
		\]
		after observing that $\ah_n$ is real and $\ah_0=0$. Furthermore, our stability condition
		requires that for all $n$, $(-\ah_n)\geq 0$, implying that the sum above is the $\ell^1$-norm
		of $\{\ah_n\}$. Putting this together with the bound~\req{inftybound} gives
		\be
			 \|a\|_\infty
			~\leq~
			 \frac{1}{M} \|\ah\|_1
			~=~
			2~ \expec{u^2_k}.
		   \label{conteffbound.eq}
		\ee

		In the vehicular formations case, the dynamics are given by~\req{vfsys} together
		with the output equation
		\[
			u ~=~ \obt{G}{F} \tbo{\xt}{\vt}.
		\]
		Formula~\req{Vv} is not applicable here since the output depends on all states,
		but the $H^2$ norm for this system can be calculated in a manner similar to the proof
		of Lemma~\ref{V.lemma}. They Lyapunov equation in this case becomes
		\beas
			\lefteqn{
			\tbt{0}{\Gh_n^*}{I}{\Fh_n^*} \tbt{\Xh_n}{\Zh_n}{\Zh_n^*}{\Yh_n} +
			\tbt{\Xh_n}{\Zh_n}{\Zh_n^*}{\Yh_n} \tbt{0}{I}{\Gh_n}{\Fh_n} 			 }
			~~~~~~~~~~~~~~~~~~~~~~~~~~~~~~~~~
			&  & 			\\
			& =  &  -\tbt{\Gh_n^*\Gh_n}{\Gh_n^* \Fh_n}{\Fh_n^* \Gh_n}{\Fh_n^* \Fh_n},
		\eeas
		from which we extract the matrix equations
		\begin{eqnarray*}
			\Gh_n^* \Zh_n^* ~+~\Zh_n \Gh_n  & = & -\Gh_n^*\Gh_n			\\
			\Gh_n^*\Yh_n ~+~\Xh_n ~+~\Zh_n \Fh_n & = & -\Gh_n^* \Fh_n 			\\
			\Zh_n ~+~\Fh_n^*\Yh_n ~+~\Zh_n^* ~+~\Yh_n \Fh_n &=& - \Fh_n^* \Fh_n  .
		\end{eqnarray*}
		Only the first and last equation need be solved since we are only interested in
		$\mbox{tr}(\Yh_n)$. All of the above are $d\times d$ diagonal matrices with equal
		entries, so we solve the equations in terms of a single entry as
		\begin{eqnarray*}
			2\gh_n \zh_n  ~ = ~ -\gh_n^2	~~~&\Rightarrow&~~~
						\zh_n ~=~ -\frac{1}{2} \gh      	\\
			2\fh_n \yh_n ~=~ -\fh_n^2 ~-~ 2\zh_n
									&\Rightarrow&
				\yh_n ~=~  \frac{1}{2} \left( -\fh_n ~+~ \frac{1}{\fh_n} \gh_n \right) .
		\end{eqnarray*}
		
		The $H^2$ norm of the system is then
		\begin{eqnarray*}
			\sum_{0\neq n\in\Z_N^d} \mbox{tr}(\Yh_n)
		& = &  \frac{d}{2}  \sum_{0\neq n\in\Z_N^d}  \left( -\fh_n ~+~ \frac{1}{\fh_n} \gh_n \right)  \\
		& = &  \frac{d}{2} \left( \|\fh\|_1 + \|\frac{1}{\fh}~\gh\|_1 \right),
		\end{eqnarray*}
		where the last equation follows from the stability conditions $\fh_n\leq 0$, $\gh_n\leq 0$.
		This inequality has two consequences after observing that
		$\sum_{n\in\Z_N^d} \mbox{tr}(\Yh_n) ~=~ M \expec{u_k^2}$ and using~\req{inftybound}
		\begin{eqnarray}
			\|f\|_\infty  & \leq &  \frac{2}{d} ~ \expec{u_k^2}
			   \label{f_bound.eq}    \\
			 \|\frac{1}{\fh}~\gh\|_1  & \leq &   \frac{2}{d}M ~ \expec{u_k^2}    .
			    \label{g/f_bound.eq}
		\end{eqnarray}
		The second inequality can be used to bound $\|g\|_\infty$ as follows.
		First
		\[
			 \|\frac{1}{\fh}~\gh\|_1  ~~\geq~~ \|\gh\|_1 ~\min_n |\frac{1}{\fh_n} |
			 ~~=~~  \|\gh\|_1 \frac{1}{\|\fh\|_\infty} .
		\]
		An upper bound on $\|\fh\|_\infty$ is derived from
		\[
			\|\fh\|_\infty ~~\leq~~ \|f\|_1  ~~\leq~~ (2q)^d \|f\|_\infty,
		\]
		where the last inequality follows from the locality assumption on $f$. Combining these
		last two bounds with~\req{g/f_bound} yields
		\[
			\frac{2}{d} M ~\expec{u_k^2}
			~\geq~ \frac{1}{(2q)^d \|f\|_\infty}  ~\|g\|_1
			~\geq~ \frac{1}{(2q)^d \|f\|_\infty}  ~M\|g\|_\infty,
		\]
		which when combined with~\req{f_bound} gives
		\[
			\|g\|_\infty ~\leq~  \frac{2(2q)^d}{d} ~\|f\|_\infty ~ \expec{u_k^2}
			~\leq~  \cB_g ~\left( \expec{u_k^2} \right)^2 .
		\]


\begin{thebibliography}{10}
\providecommand{\url}[1]{#1}
\csname url@rmstyle\endcsname
\providecommand{\newblock}{\relax}
\providecommand{\bibinfo}[2]{#2}
\providecommand\BIBentrySTDinterwordspacing{\spaceskip=0pt\relax}
\providecommand\BIBentryALTinterwordstretchfactor{4}
\providecommand\BIBentryALTinterwordspacing{\spaceskip=\fontdimen2\font plus
\BIBentryALTinterwordstretchfactor\fontdimen3\font minus
  \fontdimen4\font\relax}
\providecommand\BIBforeignlanguage[2]{{%
\expandafter\ifx\csname l@#1\endcsname\relax
\typeout{** WARNING: IEEEtran.bst: No hyphenation pattern has been}%
\typeout{** loaded for the language `#1'. Using the pattern for}%
\typeout{** the default language instead.}%
\else
\language=\csname l@#1\endcsname
\fi
#2}}

\bibitem{shldeshedtomwalzhamcmpenshemck91}
S.~Shladover, C.~Desoer, J.~Hedrick, M.~Tomizuka, J.~Walrand, W.-B. Zhang,
  D.~McMahon, H.~Peng, S.~Sheikholeslam, and N.~McKeown, ``Automated vehicle
  control developments in the path program,'' \emph{IEEE Transactions on
  Vehicular Technology}, vol.~40, no.~1, pp. 114--130, February 1991.

\bibitem{swahed99}
D.~Swaroop and J.~K. Hedrick, ``Constant spacing strategies for platooning in
  automated highway systems,'' \emph{Transactions of the ASME. Journal of
  Dynamic Systems, Measurement and Control}, vol. 121, no.~3, pp. 462--470,
  September 1999.

\bibitem{swahed96}
D.~Swaroop and J.~K. Hedrick, ``String stability of interconnected systems,''
  \emph{IEEE Transactions on Automatic Control}, vol.~41, no.~2, pp. 349--357,
  March 1996.

\bibitem{melkuo71}
S.~M. Melzer and B.~C. Kuo, ``Optimal regulation of systems described by a
  countably infinite number of objects,'' \emph{Automatica}, vol.~7, no.~3, pp.
  359--366, May 1971.

\bibitem{levath66}
W.~S. Levine and M.~Athans, ``On the optimal error regulation of a string of
  moving vehicles,'' \emph{IEEE Transactions on Automatic Control}, vol. AC-11,
  no.~3, pp. 355--361, July 1966.

\bibitem{marcorbull07}
S.~Martinez, J.~Cortes, and F.~Bullo, ``Motion coordination with distributed
  information,'' \emph{Control Systems Magazine}, vol.~27, no.~4, pp. 75--88,
  2007.

\bibitem{X07}
L.~Xiao, S.~Boyd, and S.-J. Kim, ``Distributed average consensus with
  least-mean-square deviation,'' \emph{Journal of Parallel and Distributed
  Computing}, vol.~67, pp. 33--46, 2007.

\bibitem{patbamela}
S.~Patterson, B.~Bamieh, and A.~{El Abbadi}, ``Distributed average consensus
  with stochastic communication failures,'' \emph{Submitted to IEEE Trans. Aut.
  Cont.}, 2008.

\bibitem{T84}
J.~N. Tsitsiklis, ``Problems in decentralized decision making and
  computation,'' Ph.D. dissertation, Massachusetts Institute of Technology,
  1985.

\bibitem{B90}
J.~E. Boillat, ``Load balancing and poisson equation in a graph,''
  \emph{Concurrency: Practice and Experience}, vol.~2, no.~4, pp. 289--313,
  1990.

\bibitem{jadlinmor03}
A.~Jadbabaie, J.~Lin, and A.~S. Morse, ``Coordination of groups of mobile
  autonomous agnets using neaqrst neighbor rules,'' \emph{IEEE Transactions on
  Automatic Control}, vol.~48, no.~6, pp. 988--1001, 2003.

\bibitem{PBE06}
S.~Patterson, B.~Bamieh, and A.~{El Abbadi}, ``Brief announcement: Convergence
  analysis of scalable gossip protocols,'' in \emph{20th Int. Symp. on
  Distributed Computing}, 2006, pp. 540--542.

\bibitem{carli2007cca}
R.~Carli, F.~Fagnani, A.~Speranzon, and S.~Zampieri, ``{Communication
  constraints in the average consensus problem},'' \emph{Automatica}, 2007.

\bibitem{patbam11cdc}
S.~Patterson and B.~Bamieh, ``Network coherence in fractal graphs,'' in
  \emph{Proceedings of the 50th IEEE Conference on Decision and Control and
  European Control Conference}, 2011, pp. 6445--6450.

\bibitem{barhes07}
P.~Barooah and J.~P. Hespanha, ``Estimation on graphs from relative
  measurements: Distributed algorithms and fundamental limits,'' \emph{IEEE
  Control Systems Magazine}, vol.~27, no.~4, 2007.

\bibitem{barhes08}
P.~Barooah and J.~P. Hespanha, ``Estimation from relative measurements:
  Electrical analogy and large graphs,'' \emph{IEEE Transactions on Signal
  Processing}, vol.~56, no.~6, pp. 2181--2193, 2008.

\bibitem{barhes06cdc}
P.~Barooah and J.~P. Hespanha, ``Graph effective resistance and distributed
  control: Spectral properties and applications,'' in \emph{Proc. of the IEEE
  Conf. on Decision and Control}, 2006.

\bibitem{doyle1984rwa}
P.~Doyle and J.~Snell, ``{Random Walks and Electric Networks, ser},'' \emph{The
  Carus Mathematical Monographs. Washington DC: The Mathematical Association of
  America}, 1984.

\bibitem{cserti2000alg}
J.~Cserti, ``{Application of the lattice Green's function for calculating the
  resistance of an infinite network of resistors},'' \emph{American Journal of
  Physics}, vol.~68, p. 896, 2000.

\bibitem{jovbamTAC05platoons}
M.~R. Jovanovi\'c and B.~Bamieh, ``On the ill-posedness of certain vehicular
  platoon control problems,'' \emph{IEEE Trans. Automat. Control}, vol.~50,
  no.~9, pp. 1307--1321, 2005.

\bibitem{bampagdah02}
B.~Bamieh, F.~Paganini, and M.~A. Dahleh, ``Distributed control of
  spatially-invariant systems,'' \emph{IEEE Transactions on Automatic Control},
  vol.~47, no.~7, pp. 1091--1107, July 2002.

\bibitem{cargarzam09ita}
R.~Carli, F.~Garin, and S.~Zampieri, ``Quadratic indicies for the analysis of
  consensus algorithms,'' in \emph{Proc. of the 2009 Information Theory and
  Applications Workshop}, http://ita.ucsd.edu/workshop/09/talks/, 2009.

\bibitem{barmehhes09}
P.~Barooah, P.~G. Mehta, and J.~P. Hespanha, ``Mistuning-based decentralized
  control of vehicular platoons for improved closed loop stability,'' \emph{To
  appear in IEEE Transactions on Automatic Control}, 2009.

\end{thebibliography}
\end{document}